\documentclass[twocolumn]{autart}

\usepackage[round,authoryear]{natbib}
\usepackage[normalem]{ulem} %% for crossing out texts
\usepackage{soul}  %% for crossing out texts
\usepackage{dsfont}
\usepackage[latin9]{inputenc}
\usepackage{verbatim}
\usepackage{amsmath}
\usepackage{amssymb}
\usepackage{graphicx}
\usepackage{graphics}
\usepackage{color}
\usepackage{float}
\usepackage{enumerate}
\usepackage{algorithm}
\usepackage{mathrsfs}
\usepackage{color,xspace}
\usepackage{psfrag}
\usepackage{enumitem}
\usepackage{subcaption}
\usepackage{algpseudocode}
\usepackage{lipsum}
\usepackage{cancel}

\makeatletter
\theoremstyle{plain}
\theoremstyle{remark}
\theoremstyle{definition}
\newtheorem{problem}{\protect\problemname}
\theoremstyle{definition}
\theoremstyle{plain}
\theoremstyle{plain}
\theoremstyle{plain}
\usepackage{amsfonts}
\usepackage{amsmath}
\usepackage{amssymb}
\usepackage{epsfig}
\usepackage{graphics}
\usepackage{graphicx}
\usepackage{needspace}
\usepackage{color, colortbl}
\definecolor{Gray}{gray}{0.9}

\epsfclipon
\newtheorem{theorem}{Theorem}[section]

\newtheorem{corollary}[theorem]{Corollary}
\newtheorem{criterion}{Criterion}
\newtheorem{definition}{Definition}[section]

\newtheorem{lemma}[theorem]{Lemma}

\newtheorem{proposition}{Proposition}
\newtheorem{remark}{Remark}[section]

\makeatother

\providecommand{\problemname}{Problem}

\makeatother

\providecommand{\problemname}{Problem}

\DeclareMathOperator*{\argmin}{arg\,min}
\DeclareMathOperator*{\argmax}{arg\,max}

\def\bomega{\mathbf{\boldsymbol{\omega}}}
\def\btheta{\mathbf{\boldsymbol{\theta}}}
\def\bx{\mathbf{x}}
\DeclareMathOperator{\Tr}{Tr}
\begin{document}

\begin{frontmatter}
    
    \title{Optimal network design for synchronization of coupled oscillators} % Title, preferably not more 
    
    \thanks[footnoteinfo]{This paper was not presented at any IFAC 
        meeting. Corresponding author Mahyar Fazlyab.}
    
    \author[UPENN]{Mahyar Fazlyab}\ead{mahyarfa@seas.upenn.edu},    % Add the 
    \author[ETH]{Florian D\"{o}rfler}\ead{dorfler@ethz.ch},               % e-mail address 
    \author[UPENN]{Victor M. Preciado}\ead{preciado@seas.upenn.edu}  % (ead) as shown
    
    \address[UPENN]{University of Pennsylvania, Department of Electrical and Systems Engineering, United States}  % Please supply                                              
    \address[ETH]{ETH Zurich, Automatic Control Laboratory, Switzerland}             % full addresses
    %\address[UPENN]{University of Pennsylvania, Department of Electrical and Systems Engineering, United States}

                                 % keyword list or with the 
    % help of the Automatica 
    % keyword wizard
    
\begin{abstract}
This paper studies the problem of designing networks of nonidentical coupled oscillators in order to achieve a desired level of \emph{phase cohesiveness}, defined as the maximum asymptotic phase difference across the edges of the network. 
In particular, we consider the following two design problems: \emph{(i)} the \emph{nodal-frequency design problem}, in which we tune the natural frequencies of the oscillators given the topology of the network, and (\emph{ii}) the \emph{(robust) edge-weight design problem}, in which we design the edge weights assuming that the natural frequencies are given (or belong to a given convex uncertainty set).
For both problems, we optimize an objective function of the design variables while considering a desired level of phase cohesiveness as our design constraint.
This constraint defines a convex set in the nodal-frequency design problem.
In contrast, in the edge-weight design problem, the phase cohesiveness constraint yields a non-convex set, unless the underlying network is either a tree or an arbitrary graph with identical edge weights. We then propose a convex semidefinite relaxation to approximately solve the (non-convex) edge-weight design problem for general (possibly cyclic) networks with nonidentical edge weights. 
We illustrate the applicability of our results by analyzing several network design problems of practical interest, such as power re-dispatch in power grids, sparse network design, (robust) network design for distributed wireless analog clocks, and the detection of edges leading to the Braess' paradox in power grids.
\end{abstract}

\end{frontmatter}

%%%%%%%%%%%%%%%%%%%%%%%%%%%%%%%%%%%%%%%%%%%%%%%%%%%%%%%%%%%%%%%%%%%%%%%%%%%%%%%%
%!TEX root = main_Automatica.tex
\section{Introduction}
The analysis of synchronization in networks of coupled oscillators is one of the most fundamental problems in the field of networked dynamical systems. Networks of coupled oscillators present a rich dynamic behavior, as reported in the vast literature on this topic; see, for example, \cite{FD-FB:13b} and references therein. Many complex artificial and natural systems can be modeled as networks of coupled oscillators, such as pacemaker cells in the heart, neurons in the brain, clocks in computing networks, mobile sensor networks, and power grids.
Considerable research in this field has been focused on studying the effect of network structure, coupling strengths, and nodal dynamics on the ability of a network of oscillators to synchronize \citep{di2007effects,sorrentino2007synchronizability,menck2014dead}. Various metrics have been proposed in the literature to quantify and optimize the synchronization performance. A broad class of these metrics focuses on the transient response, such as the ability of the network to resynchronize after perturbations \citep{donetti2005entangled,motter2005network,pecora1998master,kempton2015self}. In this context, synchronizability can be characterized by either the required effort to synchronize the network \citep{sjodin2014price}, the speed of convergence to the synchronization manifold \citep{xiao2004fast,fardad2014optimal}, or the range of coupling values for which a network with uniform coupling strengths would synchronize \citep{pecora1998master}. Using the master stability framework, proposed in \citep{pecora1998master}, it was shown that the Laplacian algebraic connectivity and the Laplacian eigenratio are two network-dependent measures able to capture the synchronizability of a network of \emph{identical}  coupled oscillators. Based on this connection, we find in the literature several works aiming to optimize the synchronizability of a network of identical coupled oscillators using the Laplacian matrix \citep{pecora1998master,nishikawa2006maximum,donetti2005entangled,rad2008efficient,motter2005network,motter2013spontaneous,kempton2015self,skardal2015control,fardad2014design,clark2014global,mousavi2016koopman,siami2016tractable}. 
%{\color{blue}In particular, Nishikawa et al. \citep{nishikawa2006maximum} developed an extension of the master stability framework to the case of non-diagonalizable Laplacian matrices, and found a simple condition on the eigenvalues of the Laplacian that either maximizes the ratio between the maximum and minimum stabilizing coupling strength, or minimizes the synchronization cost (defined as the sum of the coupling strengths). In \citep{donetti2005entangled,rad2008efficient}, the authors proposed a collection of rewiring algorithms to enhance synchronizability. In the context of power grids, the authors in \citep{motter2013spontaneous} considered a power system with tunable parameters and used the maximal Lyapunov exponent to characterize the local stability of the swing equation in terms of these parameters. They used this characterization to specify the parameter values in order to enhance the synchronizability of the network.} 

In \cite{FD-MC-FB:11v-pnas}, the concept of \emph{phase cohesiveness}, defined as the maximum steady-state phase difference across the edges of a network, was proposed as a synchronization metric in networks of \emph{nonidentical} coupled oscillators. This metric explicitly accounts for the simultaneous effect of the network topology, the coupling strengths, and the nodal dynamics on the local stability of the synchronous solution. This paper adopts phase cohesiveness as a synchronization measure in order to develop an optimization framework for designing the parameters of a network of coupled oscillators. As described in $\S$\ref{sub:Synchronization-in-Networks}, the oscillators in our network are modeled using the \emph{swing equation}, widely used in the analysis of power grids \citep{bergen1981structure}.
Specifically, we address the following design problems: 
\begin{enumerate}[leftmargin=*]
    \item  \emph{Design of natural frequencies}: In this problem, we assume that the network structure and the coupling strengths are given. The network designer is able to tune the natural frequencies of each oscillators by incurring a cost. The objective is to minimize the total tuning cost while guaranteeing a desired level of phase cohesiveness.
    \item \emph{Design of link weights}: In this second problem, we assume that the natural frequencies of the oscillators belong to a given polyhedral uncertainty set. The network designer is able to tune edge weights by incurring a cost. The goal is to design the edge weights while guaranteeing a desired level of phase cohesiveness for all possible realizations of the natural frequencies in the uncertainty set. 
\end{enumerate}
The framework herein proposed can be used in a wide range of practical applications, namely, prevention of cascading failures in power grids \citep{linnemann2011modeling}, optimal design of electrical infrastructure upgrades, sparsity promoting network design \citep{Siami2015246,dhingra2012identifying,lin2012identification}, and detection of links inducing the Braess' paradox (i.e., the counter-intuitive phenomenon of losing synchronization as the result of adding new edges \citep{witthaut2012braess}). We will discuss some of these applications in $\S$\ref{section:NUMERICAL SIMULATIONS}.

The rest of the paper is organized as follows. $\S$\ref{sub:Synchronization-in-Networks} provides some background on the synchronization problem. $\S$\ref{section:OPTIMAL COHESIVENESS} develops an optimization framework to solve the frequency design problem. The (robust) weight design problem is solved in $\S$\ref{subsection: weight design}. Illustrative examples are presented in $\S$\ref{section:NUMERICAL SIMULATIONS}. Concluding remarks are drawn in $\S$\ref{section:CONCLUSIONS}.

\textit{Notation}: Let $\mathbb{R}$, $\mathbb{R}_{+}$, and $\mathbb{R}_{++}$ be the set of real, nonnegative, and strictly positive numbers. Let $\mathbf{1}_n$
and $\mathbf{0}_n$ be the $n$-dimensional vectors of unit and zero entries. %$I_n$ denotes the $n$-dimensional identity matrix and $O_n$ is the $n\times n$ matrix of all zeros. 
The set $\{1,\ldots,n\}$ is denoted by $[n]$. 
%{\color{blue} Given an $n$-tuple $(x_1, . . . , x_n)$, let $\bx \in \mathbb{R}^n$ be the associated column vector.} 
The infinity norm of $\bx \in \mathbb{R}^n$ is denoted as $\|\bx\|_{\infty}=\max_i |x_i|$, the $\ell_1$ norm as $\|\bx\|_1=\sum_{i=1}^{n} |x_i|$, and the $\ell_0$ norm as $\|\bx\|_0=\mbox{card}(\{i\in[n] \colon x_{i} \neq 0 \})$, where $\mbox{card}(\cdot)$ denotes the cardinality of a set. For $\bx,\mathbf{y} \in \mathbb{R}^n$, the inequality $\bx \leq \mathbf{y}$ is component-wise. %We define the vector-valued function $\mathbf{sin}(\bx) = (\sin(x_1),\cdots,\sin(x_n))^\top$.}
%For a given matrix $A \in \mathbb{R}^{m \times n}$, the kernel of $A$ is the set $\ker(A)=\{\bx \in \mathbb{R}^n: A\bx=\mathbf{0}_m\}$. 
We denote by $\mathbb{S}^{n\times n}$ the set of $n\times n$ real, symmetric matrices. For square matrices $A$ and $B$, we write $A \succeq B$ if and only if $A-B$ is positive semidefinite.

\textit{Elements of algebraic graph theory}: A graph is defined as $G=\left(\mathcal{V},\mathcal{E}\right)$,
where $\mathcal{V}$ is a set of $n$ nodes and $\mathcal{E}$
is a set of $m$ undirected edges. We assume that the graph
is connected and has no self-loops. We consider graphs with weights
associated to both edges and nodes. We denote the weight of an edge $e=\left\{ i,j\right\} \in\mathcal{E}$
as $w_{e}=w_{ij}$. The \emph{weighted adjacency matrix} of an undirected graph $G$, denoted by $A=\left[a_{ij}\right]$,
is an $n\times n$ symmetric matrix defined entry-wise as $a_{ij}=w_{ij}$
if $\left\{ i,j\right\} \in\mathcal{E}$, and $a_{ij}=0$, otherwise.
The \emph{weighted Laplacian matrix} of $G$ is defined as $L=\mbox{diag}\left(A \mathbf{1}_n\right)-A$. For an edge $e=\{i,j\} \in \mathcal{E}$, we define $\mathbf{b}_e \in \mathbb{R}^n$ with $b_{e,i}=1,\ b_{e,j}=-1$(or $b_{e,i}=-1,\ b_{e,j}=1$) and all other entries equal to zero. The incidence matrix $B \in \mathbb{R}^{n \times m}$ is the matrix with $e$-th column $\mathbf{b}_e$. 
%{\color{blue} The (unweighted) \emph{incidence matrix} of $G$ is defined for a directed and labeled version of the undirected graph $G$, as follows. First, we label each edge in the graph with a unique label $e\in\left\{ 1,\ldots,m\right\} $ and assign an arbitrary direction to it. In other words, we substitute each undirected edge $\left\{ i,j\right\} \in\mathcal{E}$ for an ordered pair $\left(i,j\right)$, in arbitrary order. For the ordered pair $\left(i,j\right)$, we say that $i$ is the source and $j$ is the sink of the directed edge. For the resulting directed graph, the incidence matrix $B=\left[b_{ie}\right]\in\mathbb{R}^{n\times m}$ is defined component-wise as $b_{je}=1$ if $j$ is the sink node of edge $e$, $b_{ie}=-1$ if $i$ is the source node of edge $e$, and $b_{ke}=0$ otherwise. For $\bx\in\mathbb{R}^{n}$, notice that $\left(B^{\top}\bx\right)_{e}=x_{j}-x_{i}$ for a link $\left(i,j\right)$ labeled $e$.}
For a weighted graph, we define the edge-weight vector $\mathbf{w}=\left(w_{1},\cdots,w_{m}\right)^{\top}$,
where $w_{e}$ is the weight of the edge labeled $e$. The Laplacian matrix of the weighted graph can be written as $L(\mathbf{w})=B\mbox{diag}\left(\mathbf{w}\right)B^{\top}$. 
%{\color{blue} If the graph is connected, we have that $\ker\left(B^{\top}\right)=\ker\left(L\right)=\mbox{span}\left({\mathbf{1}_{n}}\right)$, and the second smallest eigenvalue of the Laplacian, called the \emph{algebraic connectivity}, is strictly positive \cite{biggs1993algebraic}.} 
The \emph{Moore-Penrose pseudoinverse} of the Laplacian is defined as $L(\mathbf{w})^{\dagger}=\left(L(\mathbf{w})+\frac{1}{n}\mathbf{1}_n\mathbf{1}_n^{\top}\right)^{-1}-{\frac{1}{n} \mathbf{1}_n\mathbf{1}_n^{\top}}$. \ifx Equivalently, if the eigenvalue decomposition of $L$ is given by $L=U\mbox{diag}(\{0,\lambda_2,...,\lambda_n\})U^\top $, its Moore-Penrose pseudo inverse is given by $L(\mathbf{w})^{\dagger}=U\mbox{diag}(\{0,\lambda_2^{-1},...,\lambda_n^{-1}\})U^\top $ \cite{biggs1993algebraic}.\fi For any connected graph with $n$ vertices, the identity $L(\mathbf{w})L(\mathbf{w})^{\dagger}=L(\mathbf{w})^{\dagger}L(\mathbf{w})=I_n-\tfrac{1}{n} \mathbf{1}_n\mathbf{1}_n^{\top}$ holds.
%Finally we denote the complement of graph $G$ with $\overline{G}$ and denote its incidence matrix as $\overline{B}$.
%!TEX root = main_Automatica.tex

%\section{Background and Preliminaries}
%\label{section:BACKGROUND AND PROBLEM STATEMENT}
%
%\subsection{\label{sub:Graph-Theory}Notation}

\section{\label{sub:Synchronization-in-Networks}Synchronization in networks of heterogeneous oscillators}
Consider a partition $\{\mathcal{V}_1,\mathcal{V}_2\}$ of the set of $n$ nodes in a connected, weighted, undirected graph ${G}(\mathcal{V},\mathcal{E})$. The state of each node $i \in \mathcal{V}$ is represented by an angular position $\theta_i \in \mathbb{R}$ whose dynamics is described by the following set of differential equations:
\begin{subequations}\label{eq:Kuramoto model}
	\begin{align}  
	m_i \ddot{\theta}_i + d_i \dot{\theta}_i &= \omega_{i} - \sum_{j=1}^{n}a_{ij}\sin\left(\theta_{i}-\theta_{j}\right),i\in\mathcal{V}_1, \\
	d_i \dot{\theta}_i &= \omega_{i} - \sum_{j=1}^{n}a_{ij}\sin\left(\theta_{i}-\theta_{j}\right), i\in\mathcal{V}_2.
	%\dot{\theta}_{i}(t)=\omega_{i} - \sum_{j=1}^{n}a_{ij}\sin\left(\theta_{i}(t)-\theta_{j}(t)\right),\forall i\in\mathcal{V}, 
	\end{align}
\end{subequations}
Here, $\mathcal{V}_1$ is a subset of oscillators following a second-order dynamics with inertia $m_i>0$ and damping coefficient $d_i>0$, and $\mathcal{V}_2$ is a subset of oscillators with a first-order dynamics; $\omega_{i} \in \mathbb{R}$ is the natural frequency of the $i$-th oscillator (which corresponds to power generation/consumption in generator/load buses), and $a_{ij}\geq0$
is the $(ij)$-th entry of the weighted adjacency matrix of ${G}(\mathcal{V},\mathcal{E})$. The dynamics in \eqref{eq:Kuramoto model} represents the swing dynamics for a structure-preserving lossless power network with constant voltage magnitudes at the buses \citep{bergen1981structure}. This dynamics can be written in matrix form as
\begin{equation} \label{eq: Kuramoto model-matrix}
M\ddot{\btheta} + D\dot{\btheta}=\mathbf{f}(\btheta) = \bomega-BW \mathbf{sin}(B^\top \btheta),
\end{equation}
where $\btheta=(\theta_1,\cdots,\theta_n)^\top $,  $\bomega=(\omega_1,\cdots,\omega_n)^\top $, $M=\mbox{diag}(\{m_i\}_{i \in \mathcal{V}_1},\mathbf{0}_{|\mathcal{V}_2|})$ is the diagonal matrix of inertias, $D=\mbox{diag}(\{d_i\}_{i\in \mathcal{V}})$ is the diagonal matrix of damping coefficients, $B$ is the incidence matrix of $G$, $W=\mbox{diag}\left(\mathbf{w}\right)$ is the diagonal matrix of edge weights, and $\mathbf{w}=(w_1,\ldots,w_m)^\top$ where $w_e>0$ is the weight of the $e$-th edge in the graph. The special case $\mathcal{V}_1=\emptyset$ and $d_i=1$ corresponds to the classical  Kuramoto model \citep{acebron2005kuramoto}. The following definition characterizes the notion of synchronization for \eqref{eq: Kuramoto model-matrix}.
%
%For future references, the model (\ref{eq:Kuramoto model}) can be expressed in vector form as follows:
%\begin{equation}
%\dot{\theta}=\omega-BW\sin(B^\top \theta)
%\end{equation}
%where $\theta=[\theta_1,\cdots,\theta_n]^\top $, $\omega=[\omega_1,\cdots,\omega_n]^\top $, and $W=\mbox{diag}\left(\mathbf{w}\right)$ is the diagonal matrix of weights.
    \begin{definition} \label{def: frequency-synchronization}
        {\rm A solution $\btheta(t)$ to the coupled oscillator model (\ref{eq:Kuramoto model}) is said to be \emph{frequency-synchronized} if $\lim\limits_{t \to \infty} |\theta_i(t)-\theta_j(t)|(\mbox{mod} \ 2\pi)=\varphi^\star_{ij}, $ for all $\{i,j\}\in \mathcal{E}$ and some $\varphi^\star_{ij} \in [0,2\pi)$. Furthermore, if $\varphi^\star_{ij}=0$ for all  $\{i,j\}\in \mathcal{E}$, the solution is said to be \emph{phase-synchronized}.}
    \end{definition}
Phase synchronization can only be achieved if all the natural frequencies are identical. In contrast, if the natural frequencies are not all identical, the network can only achieve frequency synchronization. For a frequency-synchronized solution, the angular velocities of the oscillators converge towards a common asymptotic frequency given by $ \omega_{s} = \sum_{i=1}^{n}\omega_{i}/\sum_{i=1}^{n} d_i$ \cite[$\S$ 5.2]{FD-FB:10w}. Thus, the frequency-synchronized solution satisfies $\lim\limits_{t \to \infty} (\btheta(t)-\btheta_s(t))(\mbox{mod}\ 2\pi) =\mathbf{0}_n$, where $\btheta_s(t)=(\omega_s t) \mathbf{1}_n+\btheta^{\star}$ for some $\btheta^{\star} \in \mathbb{R}^n$ such that $M\ddot{\btheta}_s + D\dot{\btheta}_s=\mathbf{f}(\btheta_s)$.  It then follows from Definition \ref{def: frequency-synchronization} that a frequency-synchronized solution $\btheta(t)$ satisfies $\lim\limits_{t \to \infty} |\theta_i(t)-\theta_j(t)| = |\theta_i^{\star}-\theta_j^{\star}|, \forall \{i,j\}\in \mathcal{E}$.
\begin{definition}
    {\rm For any frequency-synchronized solution $\btheta_s(t)=(\omega_s t) \mathbf{1}_n+\btheta^{\star}$ of \eqref{eq: Kuramoto model-matrix}, the corresponding phase cohesiveness is defined as}
    \begin{align} %\label{eq: phase_cohesiveness_definition}
    \varphi(B,\mathbf{w},\bomega)&= \underset{ \{i,j\} \in \mathcal{E}}{\max} \lim\limits_{t \to \infty} |\theta_i(t)-\theta_j(t)| (\emph{mod} \ 2\pi)\nonumber \\
    &= \|B^\top \btheta^{\star}\|_{\infty} \ (\emph{mod} \ 2\pi).
    \end{align}
\end{definition}
\vspace{-5mm}
Without loss of generality, we can assume that $\omega_s=0$ by introducing a rotational reference frame in which $\omega_s = 0$. It then follows that $\btheta_s(t)=\btheta^{\star}$ and $\mathbf{0}_n=M\ddot{\btheta}_s + D\dot{\btheta}_s=\mathbf{f}(\btheta_s)=\mathbf{f}(\btheta^{\star})$, i.e., the frequency-synchronized solution corresponds to a fixed point of (\ref{eq: Kuramoto model-matrix}),
\vspace{-2mm}
\begin{align}
\bomega-BW\mathbf{sin}(B^\top \btheta^{\star})=\mathbf{0}_n. \label{eq: fixed point}
\end{align}
The following proposition characterizes $\btheta^{\star}$, and is a generalization of the result in \cite{FD-MC-FB:11v-pnas}.
\begin{proposition} \label{lemma: fixed point}
	Define $F \in \mathbb{R}^{m\times{(m-n+1)}}$ as a matrix whose columns span the null space of $B$ (i.e., $BF=0$). Then, for any arbitrary $r \in \mathbb{R}$,  the fixed points of \eqref{eq: Kuramoto model-matrix} satisfy the following equation
	\begin{equation} \label{eq: complete solution}
	\mathbf{sin}(B^\top \btheta^{\star})=W^{r-1}B^\top (BW^{r}B^\top )^{\dagger}\bomega + W^{-1}F\mathbf{y},
	\end{equation}    
	for some vector $\mathbf{y} \in \mathbb{R}^{(m-n+1)}$ satisfying
	\begin{align} \label{eq: y vector}
	&\|W^{r-1}B^\top (BW^{r}B^\top )^{\dagger}\bomega + W^{-1}F\mathbf{y}\|_{\infty} \leq 1, \\
	&F^\top \boldsymbol{\mathbf{sin}}^{-1}\left(W^{r-1}B^\top (BW^{r}B^\top )^{\dagger}\bomega + W^{-1}F\mathbf{y}\right)=0.\nonumber
	\end{align}
\end{proposition}
\vspace{-2mm}
\begin{pf}
	It can be verified that the first summand in the right hand side of \eqref{eq: complete solution} satisfies \eqref{eq: fixed point},
	\begin{align*}
	BW W^{r-1}B^\top (BW^{r}B^\top \!)^{\dagger}\bomega \!&= \!(BW^{r}B^\top\! )(\! BW^{r}B^\top \!)^{\dagger}\!\bomega \\
	&=(I_n-\frac{1}{n}\boldsymbol{1}_n\boldsymbol{1}_n^\top )\bomega =\bomega,
	\end{align*}
	where in the second equality, we have used the fact that $L^{(r)}L^{(r)\dagger}=I_n-\frac{1}{n}\boldsymbol{1}_n\boldsymbol{1}_n^\top $ for the Laplacian matrix $L^{(r)}=BW^{r}B^\top $. The third equality follows from the assumption $\mathbf{1}_n^\top \bomega=0$. Therefore, $\mathbf{sin}(B^\top \btheta)=W^{r-1}B^\top (BW^{r}B^\top )^{\dagger}\bomega$ is a particular solution of \eqref{eq: fixed point}. The second summand in  \eqref{eq: complete solution}, $W^{-1}F\mathbf{y}$, is a homogeneous solution of \eqref{eq: fixed point}, since $BW(W^{-1}F\mathbf{y})=0$. Since $\|\boldsymbol{\sin}(B^\top \btheta^\star)\|_{\infty} \leq 1$, we have that
	\begin{align}
	\|W^{r-1}B^\top (BW^{r}B^\top )^{\dagger}\bomega + W^{-1}F\mathbf{y}\|_{\infty} \leq 1.
	\end{align}
	Furthermore, for $\btheta^{\star}$ to be realizable from \eqref{eq: complete solution}, $\mathbf{y}$ must be chosen such that $B^\top \btheta^{\star} \in \text{Im}(B^\top )$.
%	\begin{align*}
%	B^\top \btheta^{\star}&=\boldsymbol{\mathbf{sin}}^{-1}\left(W^{r-1}B^\top (BW^{r}B^\top )^{\dagger}\bomega + W^{-1}F\mathbf{y}\right) \\ &\in \text{Im}(B^\top ).
%	\end{align*}
	Or, equivalently, since $\mbox{Im}(B^\top ) \!\perp \! \ker(B)$, we must have that
	$$
	F^\top \boldsymbol{\mathbf{sin}}^{-1}\left(W^{r-1}B^\top (BW^{r}B^\top )^{\dagger}\bomega + W^{-1}F\mathbf{y}\right)=0.
	$$
	This corresponds to the geometric constraint that the sum of the phase differences along any cycle is equal to zero. The proof is complete. $\blacksquare$
	%In Appendix \ref{appendix_lemma_fixed_point}.
\end{pf}

\ifx The general result \eqref{eq: complete solution} includes special cases for different values of $r$. For instance, when $r=0$, we recover the synchronization solution obtained by a primal-dual optimization approach ( see \cite{zhou2015new} Eq. (22) ) , and $r=1$ corresponds to the synchronization condition reported in \cite{FD-MC-FB:11v-pnas}. The next lemma, proved in \cite{FD-MC-FB:11v-pnas}, provides a sufficient condition for the stability of the fixed points satisfying \eqref{eq: fixed point}.
\begin{lemma} \label{lemma: stability of fixed points}
	Any fixed point of \eqref{eq: Kuramoto model-matrix} is locally exponentially stable if $\|B^\top \btheta^{\star}\|_{\infty} < \tfrac{\pi}{2}$.
\end{lemma}
It follows directly from the last lemma that the condition $\varphi(B,\mathbf{w},\bomega) < \pi/2$ is sufficient for local exponential stability of the fixed point. As a special case, in trees with uniform coupling weights, $F\mathbf{y}=\mathbf{0}_{m}$ and $\mathbf{w}=k\mathbf{1}_{n-1}$ where $k>0$ is the coupling strength, the phase cohesiveness reduces to $\sin(\varphi(B,k\mathbf{1}_{n-1},\bomega))=k^{-1}{B^T(BB^T) ^{\dagger}\bomega}$. Consequently, the stability condition in Lemma \ref{lemma: stability of fixed points} translates into the condition $k>\|B^T(BB^T) ^{\dagger}\bomega\|_{\infty}$ for tree graphs.\fi

It was shown in \cite{taylor2015finding} that finding the nonzero stable fixed points of \eqref{eq: Kuramoto model-matrix} over the full space of phase angles $[0,2\pi)^n$ is NP-hard. Alternatively, the following synchronization criterion, proposed in \cite{FD-MC-FB:11v-pnas}, can be used to find an upper bound for the phase cohesiveness.
\begin{criterion} \label{sync criterion} The oscillator model \eqref{eq:Kuramoto model} has a unique and stable frequency-synchronized solution $\btheta^{\star}$ such that $|\theta_i^{\star}-\theta_j^{\star}| \leq \gamma < \pi/2$ for every $\{i,j\} \in \mathcal{E}$ if
    \begin{equation} 
    \left\Vert B^{\top}L(\mathbf{w})^{\dagger}\bomega\right\Vert _{\infty}\leq\sin\left(\gamma\right).\label{eq:sync condition}
    \end{equation}%where $L=B\emph{diag}(\mathbf{w})B^\top$.
\end{criterion}
\vspace{-3mm}
The above criterion implies that when $\|B^\top L(\mathbf{w})^{\dagger} \bomega\|_{\infty}<1$, the phase cohesiveness satisfies
\begin{align} \label{eq: criterion equivalent}
\quad \varphi(B,\mathbf{w},\bomega) \leq \sin^{-1}(\| B^{\top}L(\mathbf{w})^{\dagger}\bomega\| _{\infty}).
\end{align}
In other words, \eqref{eq: criterion equivalent} provides an upper bound on the phase cohesiveness in terms of $B$, $\mathbf{w}$, and $\bomega$. In particular, the condition $\|B^\top L(\mathbf{w})^{\dagger}\bomega\|_{\infty}<1$ (corresponding to $\gamma=\pi/2$ in \eqref{eq:sync condition}) implies that $\varphi(B,\mathbf{w},\bomega) \leq \pi/2$, which guarantees local exponential stability and is a common security index for the stability of power systems \citep{FD-MC-FB:11v-pnas}. It can also be shown that the upper bound in \eqref{eq:sync condition} is tight for various topologies including trees and complete graphs. In particular, for tree graphs, we have $F=0$ in \eqref{eq: complete solution}, and by setting $r=1$, we obtain $\sin(\varphi(B,\mathbf{w},\bomega))=\| B^{\top}L(\mathbf{w})^{\dagger}\bomega\| _{\infty}$; see \cite{FD-MC-FB:11v-pnas} for further details.

\ifx
\begin{remark}{\rm
    It was shown in \cite{FD-MC-FB:11v-pnas} that the criterion \ref{sync criterion} is provably correct for various
    network topologies and weights including the extremal cases of the
    sparsest (acyclic) and densest (complete) graphs.
    By statistical means, it has also been shown that the inequality (\ref{eq:sync condition})
    is extremely accurate for a broad set of random network topologies
    and weights, as well as for various standard power network test cases
    \cite{FD-MC-FB:11v-pnas}.}
\end{remark}
\fi
\ifx
\begin{remark} \label{cor: tree graphs}
    For tree graphs, the null space of $B$ is trivial, i.e., $F=\mathbf{0}_{m \times (m-n+1)}$; hence, according to \eqref{eq: complete solution} the fixed point solution becomes independent of $r$. More specifically, using $r=1$, we obtain $\mathbf{sin}(B^\top\btheta^{\star})=B^\top L(\mathbf{w})^{\dagger}\bomega$ and therefore, the phase cohesiveness bound in \eqref{eq: criterion equivalent} is tight for acyclic graphs. On the other hand, for $r=0$, we get $\mathbf{sin}(B^\top\btheta^{\star})=W^{-1}B^\top (BB^\top)^{\dagger}\bomega$
\end{remark}
\fi
The next section uses the upper bound in \eqref{eq: criterion equivalent} to develop an optimization framework for designing the natural frequencies of the network with given edge weights.

%Therefore, we must resort to approximation methods to compute the fixed points.
%For example, it can be shown that the solution to the following quadratic optimization problem:
%\begin{align}  \label{eq: approximate solution}
%\bx^{\star}=\arg \underset{\boldsymbol{x}}{\min} & \quad \dfrac{1}{2} \boldsymbol{x}^\top W \boldsymbol{x} \nonumber \\
%\mbox{s.t. } & \quad BW\boldsymbol{x}=\bomega.
%\end{align}
%satisfies $\bx^{\star}=B^\top L(\mathbf{w})^{\dagger}\bomega$.
%The dynamics of the Kuramoto model can be interpreted as a distributed gradient descent on the total energy  \textcolor{red}{ defined how???} of the of the system \cite{FD-MC-FB:11v-pnas}. On the other hand, the optimization problem \eqref{eq: approximate solution} is indeed finding the minimum energy solution among all the fixed point equations satisfying \eqref{eq: Kuramoto model-matrix}.

\ifx
\subsection{\label{sub:sync-metric}Synchronization metrics}
According to Criterion \ref{criterion}, the infinity norm of the approximate solution \eqref{eq: approximate solution} provides an upper bound for the cohesiveness of the synchronized phases. This upper-bound is tight for most real world oscillators. Therefore, we construct the metrics of synchronization based on this upper bound.
\medskip

\begin{enumerate}[(a)]

\item \textbf{Phase cohesiveness}: quantifies an upper-bound for the maximum phase difference across the edges. The infinity norm of \eqref{eq: sin(phistar)} yields the phase cohesiveness $\Phi: \mathbb{R}^m_{+} \times \mathbb{R}^n \to \mathbb{R}_{+}$ as follows: 
    \begin{align} \label{eq: phase cohesiveness}
    \Phi(\mathbf{w},\bomega) &=\|B^\top (BWB^\top )^{\dagger}\bomega\|_{\infty},
    \end{align}
    The following definition is the consequence of relating lemma \ref{lemma: stability of fixed points} to Criterion \ref{criterion}.
\begin{definition}
    A network of oscillators with incidence matrix $B$, weight vector $\mathbf{w}$ and vector of natural frequency $\bomega$ is synchronizable if $\Phi(\mathbf{w},\bomega) < 1$.
\end{definition}
\begin{remark}
    Proof of Lemma \ref{lemma: stability of fixed points} implies that as $\Phi(\mathbf{w},\bomega) \to 1$, the second largest eigenvalue of the Jacobian in increases toward the origin, implying that phase cohesiveness captures the margin of stability for the synchronized solution.
\end{remark}

\begin{remark}
    For particular graph structures, the critical coupling \eqref{eq: Critical Coupling} can be simplified. In particular, for chain and star graphs, the critical coupling can be rewritten as
    \begin{align}
    &\Phi^{chain}(\mathbf{w},\bomega)=\underset{1 \leq j \leq n-1}{\max}\left|\sum\limits_{i=1}^{j}\dfrac{(\omega_i-\omega_{s})}{w_i}\right|, \\  &\Phi^{star}(\mathbf{w},\bomega)=\underset{2 \leq j \leq n}{\max}\left|\dfrac{\omega_i-\omega_{s}}{w_i}\right|, \nonumber
    \end{align}
    Dekker et al.\cite{dekker2013synchronization} have studied the approximation of the expected critical coupling for various classes of tree when frequencies are uniformly distributed on the unit interval.
\end{remark}

\item \textbf{Phase efficiency} quantifies the collective  performance of the whole network from an energy perspective. The phase efficiency is quantified by the weighted sum of squares of the phase differences:
\begin{align} \label{eq: phase efficiency}
\mathcal{P}(\mathbf{w},\bomega)&=\mathbf{sin}^\top (\boldsymbol{\varphi}^{\star}) Q\mathbf{sin}(\boldsymbol{\varphi}^{\star}) \nonumber \\
&=\|B^\top (BWB^\top )^{\dagger}\bomega\|^{2}_{2,Q}
\end{align}
\begin{rem}
    In power networks, it can be shown that the instantaneous resistive power loss is proportional to the square of phase differences weighted by the conductances of the transmission lines \cite{sjodin2014price}.
\end{rem}

\ifx
\begin{rem}
    For the case of uniform natural frequencies, we have that $\bomega=\alpha \mathbf{1}_n$ and since $(BW^rB^\top )^{\dagger}\mathbf{1}_n=0$, it follows from \eqref{eq: phase cohesiveness}, \eqref{eq: phase efficiency}  all performance measures approach zero as $\bomega \to \alpha 
    \mathbf{1}_n$, meaning, the network perfectly synchronizes. In this case, the coupling weights only determine the speed of synchronization.
\end{rem}    
\fi

\end{enumerate}

\fi
%!TEX root = main_Automatica.tex

\section{Cohesiveness-constrained frequency design}
\label{section:OPTIMAL COHESIVENESS}

Consider the model \eqref{eq: Kuramoto model-matrix} where $B$ and $\mathbf{w}$ are given. In the frequency design problem, our objective is to design $\bomega$, within a convex feasible set, such that a desired level of phase cohesiveness $\gamma_d \in [0, \pi/2)$ is guaranteed at a minimum design cost. This problem can be mathematically stated as follows.
 \begin{problem}\textbf{\emph{(Frequency design)}}{\rm
    \label{prob:Optimal Budget Problem-Frequency design} Assume we are given the following elements: (i) a connected undirected network with incidence matrix $B$, (ii)
    a nonnegative vector of link strengths $\mathbf{w}_{0}\in\mathbb{R}_{+}^{m}$,
    (iii) a convex frequency-tuning cost function $g_{\mathcal{V}}\left(\mathbf{\bomega}\right): \mathbb{R}^n \to \mathbb{R}$, (iv) a closed convex feasible design set $F_{\mathbf{\bomega}}\subset\mathbb{R}^{n}$, (v) a desired level of phase cohesiveness $\gamma_d \in [0,\pi/2)$, and (vi) a synchronizing frequency $\omega_s \in \mathbb{R}$. Find an optimal vector of natural frequencies, denoted by $\bomega^{\star}$, that solves
    \begin{align} \label{eq: Cohesiveness-Constrained Frequency Design}
    \bomega^{\star} \in &\underset{\bomega\in F_{\mathbf{\bomega}}}{\argmin} \  g_{\mathcal{V}}\left(\mathbf{\bomega}\right)\\
    &\mbox{s.t.}  \  \|B^\top L(\mathbf{w}_0)^{\dagger}\bomega\|_{\infty}\leq \sin(\gamma_d), \nonumber \
    \frac{1}{n}\boldsymbol{1}_n^\top   \bomega=\omega_s, \nonumber
     \end{align}
    where $L(\mathbf{w}_0)=B\mbox{diag}(\mathbf{w}_0)B^\top $. Then, by Criterion \ref{sync criterion}, the phase cohesiveness of the resulting network would satisfy 
    $\varphi(B,\mathbf{w}_0,\bomega^{\star}) \leq \gamma_d$ at a minimum cost. }
\end{problem}

Note that Problem \ref{prob:Optimal Budget Problem-Frequency design} is a convex optimization problem, since the function $\|B^\top L(\mathbf{w}_0)^{\dagger}\bomega\|_{\infty}$ can be written as the point-wise maximum of linear functions of $\bomega$; hence, the feasible set $\{\bomega \in F_{\bomega} \colon \ \|B^{\top}L(\mathbf{w}_0)^{\dagger}\bomega\|_{\infty} \leq \sin(\gamma_d) \}$ renders a convex set.

%Notice that if $\omega_s \mathbf{1}_n \in F_{\bomega}$, then $\bomega^\star = \omega_s \mathbf{1}_n$ is the trivial optimal solution. Therefore, we assume tha the feasible design set does not contain $\bomega = \omega_s \mathbf{1}_n$

%{\color{blue} Notice that in the context of power systems, the thermal limit constraints on the transmission lines are precisely equivalent to bounds on the phase cohesiveness \citep{FD-MC-FB:11v-pnas}. To respect these limits, remedial actions such as redispatch or load shedding (tuning $\bomega$) can be used to satisfy the phase cohesiveness constraint in \eqref{eq: Cohesiveness-Constrained Frequency Design} \citep{shi2015decentralized}. In this case, we can choose the cost function $g_{\mathcal{V}}(\bomega)$ to be a penalty on the amount of load (or power generation) that has to be shed (or redispatched). See Section \ref{subsection:Power Redispatch} for an application.}

\ifx
\begin{remark}
    A closely related problem to cohesiveness-constrained frequency design is the so-called budget-constrained frequency design. In this problem, the goal is to minimize the phase cohesiveness given a certain tuning budget $C$ that the network designer can invest on tuning the natural frequencies.
\end{remark}

In the next section, we address the \emph{weight design problem} where the objective is to design the coupling strengths of a network in order to achieve a desired level of phase cohesiveness at a minimum design cost. \fi
\ifx
A special case of frequency design problem is uniform weights on the edges, i.e. $\mathbf{w}=k \mathbf{1}_n$. In this case, the phase cohesiveness is simply $\Phi(k,\bomega)=k^{-1}\|B^\top (BB^\top )^{\dagger}\bomega\|_{\infty}$. The following problem addresses this case.
\begin{problem}[Joint Frequency-Weight Design for Homogeneous Coupling Weights] \label{prob:optimal critical coupling-Frequency}
    Given (\emph{i}) a graph with incidence matrix $B$ and uniform coupling strengths $\mathbf{w}=k \mathbf{1}_m$ (\emph{ii}) a frequency tuning function $g_{\mathcal{V}}\left(\bomega\right): \mathbb{R}^n \to \mathbb{R}^+$,    (\emph{iii}) a feasible design set $F_{\mathbf{\omega}}\subset\mathbb{R}^{n}$, (\emph{iv}) a desired level of phase cohesiveness ${\Phi}^{\star}\in\left[0,1\right)$, and (\emph{v}) a synchronized frequency $\omega_{s}$; find the optimal vector of natural frequencies, denoted by $\bomega^{\star}$, such that the network synchronizes at $\omega_{s}$ with cohesiveness level $\Phi^{\star}$ with minimum possible coupling. 
    
    Problem \ref{prob:optimal critical coupling-Frequency} can be formulated as follows:
    \begin{equation}
    \begin{aligned}\underset{k> 0, \bomega\in F_{\mathbf{\bomega}}}{\min} & \quad k\\
    \mbox{s.t. } & \quad \|B^\top (BB^\top )^{\dagger}\bomega\|_{\infty} \leq k \Phi^{\star}\\
    & \quad g_{\mathcal{V}}\left(\bomega\right)\leq C,\\
    %& \quad \bomega\in F_{{\omega}},\\
    \end{aligned}
    \end{equation}
\end{problem}

\begin{rem}
    As was shown in section \ref{subsection: stability analysis}, at the critical coupling, the phase cohesiveness is equal to one. Hence, Problem \ref{prob:optimal critical coupling-Frequency} with $\Phi^{\star}=1$ finds the set of feasible natural frequencies that minimize the critical coupling.
\end{rem}
\fi

%!TEX root = main_Automatica.tex

\section{\label{subsection: weight design}Cohesiveness-constrained weight design}

Consider the model \eqref{eq: Kuramoto model-matrix} with a given $B$ and $\bomega$. In the weight design problem, our objective is to design $\mathbf{w}$, within a convex feasible set $F_{\mathbf{w}} \subset \mathbb{R}^m_{+}$, such that a desired level of phase cohesiveness $\gamma_d \in [0, \pi/2)$ is achieved at a minimum cost. We assume that $\lambda_2(L(\mathbf{w}))>0$ for all $\mathbf{w}\in F_{\mathbf{w}}$, which implies that the design set excludes disconnected graphs. In most practical settings, however, the natural frequencies are uncertain. For example, in the context of power systems, the natural frequencies correspond to net power injected/consumed at the buses; thus, these values are subject to uncertainties depending on demand generation patterns. %We find another example in the context of distributed synchronization of coupled analog clocks \citep{simeone2008distributed}. The natural frequencies of the clocks are inherently uncertain due to manufacturing defects. It is therefore desirable to have a communication network capable of synchronizing the clocks to a certain level of phase cohesiveness despite the uncertainties in the natural frequencies. 
Consequently, it is of practical relevance to extend the weight design problem in order to support uncertainties in the natural frequencies. In this direction, we assume that $\bomega$ in the following polyhedral uncertainty set,
\begin{align} \label{eq: freqeuency_uncertainty_set}
\Omega=\{\bomega \in \mathbb{R}^n \colon \ C \bomega \leq \mathbf{d}\},
%%\Omega=\{\bomega \in \mathbb{R}^n \colon \ \underline{\bomega} \leq \bomega \leq \overline{\bomega}, \ \mathbf{1}_n^\top \bomega = 0\},
\end{align}
where $C \in \mathbb{R}^{p\times n}$, and $\mathbf{d} \in \mathbb{R}^p$ are given. We formalize the (robust) weight design problem next.
\begin{problem}\textbf{\emph{(Weight design)}} \label{problem: Budget-Constrained Weight Design}{\rm
    \label{prob:Optimal Cohesiveness Problem-Weight design} Assume we are given the following elements: (i) a connected undirected network with incidence matrix $B$, (ii) a polyhedral uncertainty set $\Omega \subset \mathbb{R}^n$ of possible values for the vector of natural frequencies (see \eqref{eq: freqeuency_uncertainty_set}), (iii) a convex objective function $f_{\mathcal{E}}\left(\mathbf{w}\right):\mathbb{R}_+^m \to \mathbb{R}$,
    (iv) a closed convex feasible design set $F_{\mathbf{w}}\subset\mathbb{R}_{+}^{m}$,
    and (v) a desired level of phase cohesiveness $\gamma_d \in [0,\pi/2)$. Find an optimal vector of link weights, denoted by $\mathbf{w}^{\star}$, that solves
    %
%    \begin{align}\label{eq: robust weight design}
%    	\underset{\mathbf{w}\in F_{\mathbf{w}}}{\mbox{minimize}} \ f_{\mathcal{E}}\left(\mathbf{w}\right) \quad \mbox{s.t. } \; \underset{\bomega \in \Omega}{\max} \ \varphi(B,\mathbf{w},\bomega) \leq \gamma_d.
%    \end{align}}
\begin{align}\label{eq: robust weight design}
&\underset{\mathbf{w}\in F_{\mathbf{w}}}{\mbox{minimize}} \ f_{\mathcal{E}}\left(\mathbf{w}\right) \\
& \mbox{s.t. } \underset{\bomega \in \Omega}{\max} \ \|B^\top  L(\mathbf{w})^{\dagger}\bomega \|_{\infty}\leq \sin(\gamma_d). \nonumber
\end{align}
Then, by Criterion \ref{sync criterion}, the phase cohesiveness satisfies $\varphi(B,\mathbf{w}^{\star},\bomega) \!$ $\leq \! \gamma_d$ for all $\bomega  \in \Omega$.}
\end{problem}     
    \vspace{-2mm}
    It turns out that any symmetric\footnote{A symmetric function is invariant under any permutation of its arguments \citep{borwein2010convex}.}, closed, convex function of the nontrivial eigenvalues of the Laplacian, denoted by $0<\lambda_2(\mathbf{w})\leq \ldots \leq \lambda_n(\mathbf{w})$ is a convex function of the edge weights \citep{borwein2010convex}. In particular, the algebraic connectivity $\lambda_2(\mathbf{w})$, the Laplacian eigenratio $\lambda_n(\mathbf{w})/\lambda_2(\mathbf{w})$, and the total effective resistance $R=n^{-1}\sum_{i=2}^{n} \lambda_i^{-1}(\mathbf{w})$ \citep{boyd2006convex} are tractable objective functions in our framework. The next proposition shows that the robust optimization problem \eqref{eq: robust weight design} is convex for acyclic connected networks.
    
    \begin{proposition} \label{prop: weight_design_tree}
    	{\normalfont \textbf{(Robust weight design for tree networks)}}The robust weight optimization in \eqref{eq: robust weight design} for acyclic connected topologies is equivalent to the following convex problem,
    	\begin{subequations} \label{eq: weight_design_tree}
    	\begin{align} \label{eq: weight_design_tree 1}
    	\underset{\mathbf{w}\in F_{\mathbf{w}}}{\emph{minimize}} \ f_{\mathcal{E}}\left(\mathbf{w}\right) \quad \emph{s.t.} \; \mathbf{w} \geq \underline{\mathbf{w}} \,  \sin(\gamma_d)^{-1},
    	\end{align}
    	where each component of $\underline{\mathbf{w}}=[\underline{w}_{1},...,\underline{w}_{n-1}]^\top \in \mathbb{R}_{++}^{n-1}$ is the solution to the following linear program (LP),
    	\begin{align}  \label{eq: weight_design_tree 2}
    	\underline{w}_{e} = \max_{\bomega \in \Omega} |\mathbf{u}_e^\top B^\top(BB^\top)^\dagger \bomega|, \ e\in[n-1],
    	\end{align}
     	 \end{subequations}
    	where $\mathbf{u}_e \in \mathbb{R}^m$ is the $e$-th standard unit basis vector. 
    \end{proposition}
    \begin{pf} 
    	%See Appendix \ref{appendix: tree_robust_problem}.    	
    	For acyclic topologies, the incidence matrix $B \in \mathbb{R}^{n\times (n-1)}$ is full column rank, implying that  $F=\emptyset$ in \eqref{eq: complete solution}. By setting $r=0$ and $r=1$ in \eqref{eq: complete solution}, we obtain the identity $\|B^\top  L(\mathbf{w})^{\dagger}\bomega \|_{\infty} = \|W^{-1}B^\top (BB^\top)^\dagger \bomega\|_{\infty}$. Expanding the right-hand side term yields
    	%\begin{align*} %\label{eq: phase cohesiveness tree}
    	%\sin (\varphi(B,\mathbf{w},\bomega))&=\|W^{-1}B^\top (BB^\top)^\dagger \bomega\|_{\infty},
    	%%&=\underset{e \in [m]}{\max} \frac{|\mathbf{1}_e^\top B^\top(BB^\top)^\dagger \bomega|}{w_e},
    	%\end{align*}
    	%which can be rewritten as
    	%The phase cohesiveness expression \eqref{eq: phase cohesiveness tree} can be rewritten as
    	\begin{align*}
    	\|B^\top  L(\mathbf{w})^{\dagger}\bomega \|_{\infty} =\underset{e \in [n-1]}{\max} \frac{|\mathbf{u}_e^\top B^\top(BB^\top)^\dagger \bomega|}{w_e}.
    	\end{align*}
    	Therefore, $\|B^\top  L(\mathbf{w})^{\dagger}\bomega \|_{\infty} $ is the point-wise maximum of the convex functions $w_e \mapsto |\mathbf{u}_e^\top B^\top(BB^\top)^\dagger \bomega|/w_e$ for $e\in[n-1]$. Therefore, it is convex in $\mathbb{R}^m_{++}$. By substituting this expression in the constraint of \eqref{eq: robust weight design}, interchanging the max operations, and defining $\underline{w}_{e} = \max_{\bomega \in \Omega} |\mathbf{u}_e^\top B^\top(BB^\top)^\dagger \bomega|$, we will arrive at \eqref{eq: weight_design_tree}. The proof is complete. $\blacksquare$
    	%\begin{align*} %\label{eq: robust weight design_tree}
    	%&\underset{\mathbf{w}\in F_{\mathbf{w}}}{\mbox{minimize}} \ f_{\mathcal{E}}\left(\mathbf{w}\right) \\
    	%& \mbox{s.t. }  \max_{\bomega \in \Omega} \underset{e \in [m]}{\max} \frac{|\mathbf{u}_e^\top B^\top(BB^\top)^\dagger \bomega|}{w_e},\leq \sin(\gamma_d). \nonumber
    	%\end{align*}
    	%By interchanging the max operations in the constraint, we get
    	%\begin{align*} %\label{eq: robust weight design_tree}
    	%\mathbf{w}^{\star} \in &\arg\underset{\mathbf{w}\in F_{\mathbf{w}}}{\min} \ f_{\mathcal{E}}\left(\mathbf{w}\right) \\
    	%& \mbox{s.t. } \max_{\bomega \in \Omega} \frac{|\mathbf{u}_e^\top B^\top(BB^\top)^\dagger \bomega|}{w_e},\leq \sin(\gamma_d), \forall e\in[m]. \nonumber
    	%\end{align*}
    	%By defining $\underline{w}_{e} = \max_{\bomega \in \Omega} |\mathbf{u}_e^\top B^\top(BB^\top)^\dagger \bomega|$, the proof is complete. 
    \end{pf}
    \vspace{-5mm}
    
According to Proposition \ref{prop: weight_design_tree}, for acyclic graphs, the phase cohesiveness constraint in \eqref{eq: robust weight design} translates into lower bounds on the weights. These lower bounds are the solutions of $(n-1)$ LP problems and, hence, can be computed efficiently. 

\begin{remark}
	\textbf{(Critical coupling for tree networks)} It follows from Proposition \ref{prop: weight_design_tree} that in the case of tree graphs with uniform weights (i.e., when $\mathbf{w}=w \mathbf{1}_n$, $w>0$), the inequality constraints in \eqref{eq: weight_design_tree 1} simplifies to
	\begin{equation} \label{eq: critical coupling}
	w \geq \|{B^\top(BB^\top) ^{\dagger}\boldsymbol{\omega}}\|_{\infty},
	\end{equation}
	after the substitutions $\Omega=\{\omega\}$ and $\gamma_d=\pi/2$. Since $\gamma_d=\pi/2$ is the stability threshold for synchronization (see the discussion after Criterion \ref{sync criterion}), the lower bound in \eqref{eq: critical coupling} corresponds to the minimum edge weight for which the network synchronizes, also known as the \emph{critical coupling}; see \cite{jadbabaie2004stability} and \cite{dekker2013synchronization}.
\end{remark}

%For general graphs, we use Criterion \ref{sync criterion} to replace the phase cohesiveness function in \eqref{eq: robust weight design} by its upper bound in \eqref{eq:sync condition}, resulting in the following relaxed formulation:
%    %
%    \begin{align}\label{eq: robust weight design_2}
%	&\underset{\mathbf{w}\in F_{\mathbf{w}}}{\mbox{minimize}} \ f_{\mathcal{E}}\left(\mathbf{w}\right) \\
%    & \mbox{s.t. } \underset{\bomega \in \Omega}{\max} \ \|B^\top  L(\mathbf{w})^{\dagger}\bomega \|_{\infty}\leq \sin(\gamma_d). \nonumber
%    \end{align}
%    Then, by Criterion \ref{sync criterion}, we have that $\max_{\bomega \in \Omega} \!\varphi(B,\mathbf{w}^{\star},\bomega) \!$ $\leq \! \gamma_d$, where $\mathbf{w}^\star$ is any optimal solution to \eqref{eq: robust weight design_2}.
%    
%    \ifx
%    It can be verified that $\|B^\top ( B\mbox{diag}({\mathbf{w}})B^\top )^{\dagger}\bomega \|_{\infty}$, appearing in the constraint of \eqref{eq: robust weight design}, remains unchanged by the substitution $\bomega \leftarrow \bomega-\frac{1}{n} \mathbf{1}_n \mathbf{1}_n^\top \bomega$; hence, we can include the assumption $\bomega^\top \mathbf{1}_n=0$ in the uncertainty set $\Omega$ without loss of generality.\fi

For tree graphs with identical edge weights, the constraint in the problem in \eqref{eq: robust weight design} simplifies to $\{w>0 \colon \! \max_{\bomega \in \Omega} \|B^\top  (BB^\top )^{\dagger}\bomega \|_{\infty}\leq w \sin(\gamma_d)\}$, which is a convex set; hence, the optimization problem in  \eqref{eq: robust weight design} is convex in $w$. However, for the case of trees with nonidentical edge weights, the weight design problem in \eqref{eq: robust weight design} is non-convex and typically intractable for general uncertainty sets. We show below how to use robust optimization tools \citep{bertsimas2011theory} to convert the problem into a tractable form when $\Omega$ is a polyhedron (as defined in \eqref{eq: freqeuency_uncertainty_set}). The main idea is to use duality theory to replace the subproblem $\max_{\bomega \in \Omega} \|B^\top L(\mathbf{w})^{\dagger} \bomega\|$ by its dual function, which provides a tight upper bound not involving the uncertain parameter $\bomega$. The next theorem provides a tractable formulation akin to \eqref{eq: robust weight design}.

\begin{theorem}[Robust weight design]\label{thm: robust weight design equivalent}
    Consider the optimization problem \eqref{eq: robust weight design} with $\Omega$ defined in \eqref{eq: freqeuency_uncertainty_set}. Then, the following optimization problem is equivalent to \eqref{eq: robust weight design},
    \begin{align} \label{eq: robust weight design equivalent}
    \text{\emph{min}}\ &f_{\mathcal{E}}(\mathbf{w}) \\ 
    \mbox{\emph{s.t. }}\ &\forall e \in [m], \nonumber \\
    &C^\top{\boldsymbol{\lambda}}_e \!- \!L(\mathbf{w})^{\dagger}\mathbf{b}_e=\mathbf{0}_n, \ C^\top{\boldsymbol{\gamma}}_e\!+\! L(\mathbf{w})^{\dagger}\mathbf{b}_e =\mathbf{0}_n, \nonumber \\
    &\boldsymbol{\lambda}_e^\top \mathbf{d} \leq \sin(\gamma_d), \ \boldsymbol{\gamma}_e^\top \mathbf{d}  \leq \sin(\gamma_d), \ \boldsymbol{\lambda}_e, \  {\boldsymbol{\gamma}}_e \geq \mathbf{0}_n,\nonumber
    \end{align}
    with variables $\mathbf{w}\in F_{\mathbf{w}}$, and $\boldsymbol{\lambda}_e, {\boldsymbol{\gamma}}_e \in \mathbb{R}^p_{+}, e\in[m]$. 
\end{theorem}
\vspace{-5mm}
    \begin{pf}
        %In Appendix \ref{appendix thm: robust weight design equivalent proof}.
        %\vspace{-15mm}
        We use $\|\bx\|_{\infty}=\max_{e \in [m]} \{x_e,-x_e\}, \ \bx \in \mathbb{R}^m$ to expand the subproblem $\max_{\bomega \in \Omega} \|B^\top L(\mathbf{w})^{\dagger} \bomega\|_{\infty}$ $\leq \sin(\gamma_d)$  as follows,                        
        \ifx
        \begin{align}
        \boldsymbol{\varphi}(\mathbf{w},\bomega) = B^\top  ( B\mbox{diag}({\mathbf{w}})B^\top )^{\dagger}\bomega
        \end{align}
        \begin{align*}
        &\underset{\bomega \in \Omega}{\max} \ \|\boldsymbol{\varphi}(\mathbf{w},\bomega)\|_{\infty} = \\ &\underset{\bomega \in \Omega}{\max} \quad \underset{e \in [m]}{\max} |\boldsymbol{\varphi}_e(\mathbf{w},\bomega)|= \\ 
        & \underset{\bomega \in \Omega}{\max} \quad \underset{e \in [m]}{\max} \quad {\max} \{\boldsymbol{\varphi}_e(\mathbf{w},\bomega),-\boldsymbol{\varphi}_e(\mathbf{w},\bomega)\} = \\        
        &\underset{e \in [m]}{\max} \quad \max \{ \underset{\bomega \in \Omega}{\max} \ \boldsymbol{\varphi}_e(\mathbf{w},\bomega),\underset{\bomega \in \Omega}{\max} \ -\boldsymbol{\varphi}_e(\mathbf{w},\bomega)\}
        \end{align*}
        In other words, the constraint $\underset{\bomega \in \Omega}{\max} \ \|\boldsymbol{\varphi}(\mathbf{w},\bomega)\|_{\infty} \leq t$ is equivalent to
        $\max \{ \underset{\bomega \in \Omega}{\max} \ \boldsymbol{\varphi}_e(\mathbf{w},\bomega),\underset{\bomega \in \Omega}{\max} \ -\boldsymbol{\varphi}_e(\mathbf{w},\bomega)\}$ for all $e \in [m]$. 
        \fi
        \begin{subequations}
	        \begin{align} %\label{eq: robust_weight_design_equiv_slack}
        	&  \underset{\bomega \in \Omega}{\mbox{max}}  \quad \ \, \mathbf{b}_e^\top L(\mathbf{w})^{\dagger}\bomega \leq \sin(\gamma_d), \ e\in [m], \label{eq: max constraint} \\ 
        	&  \underset{\bomega \in \Omega}{\mbox{max}} \ -\mathbf{b}_e^\top L(\mathbf{w})^{\dagger}\bomega \leq \sin(\gamma_d), \ e\in [m]. \label{eq: max constraint 2} 
        	\end{align}
        \end{subequations}
       For each fixed $\mathbf{w} \in F_{\mathbf{w}}$ and $e \in [m]$, the left-hand side of \eqref{eq: max constraint} is an LP in $\bomega \in \Omega$. For the specific choice of $\Omega$ as in \eqref{eq: freqeuency_uncertainty_set}, the $e$-th Lagrangian function is                    
        \begin{align} \label{eq: lagrangian function}
        \mathcal{L}_e(\bomega,{\boldsymbol{\lambda}}_e)&=\mathbf{b}_e^\top L(\mathbf{w})^{\dagger}\bomega+{\boldsymbol{\lambda}}_e^\top(\mathbf{d}-C\bomega),
        \end{align}            
        where ${\boldsymbol{\lambda}}_e, \in \mathbb{R}^{p}_{+}$ are the Lagrange multipliers. Notice that although the constraint set $\Omega$ is identical for all $e\in [m]$, the individual objective functions ($\mathbf{b}_e^\top L(\mathbf{w})^{\dagger}\bomega$) depend on $e$; so do the Lagrange multipliers. The dual function of \eqref{eq: lagrangian function} is defined as
        %\begin{align} \label{eq: dual function}
        $G_e(\boldsymbol{\lambda}_e)=\max_{\bomega \in \mathbb{R}^n}  \mathcal{L}_e(\bomega,{\boldsymbol{\lambda}}_e)$.
        %\end{align}
        %
        The Lagrangian function is an affine function of $\bomega$ and, hence, the dual function is infinite unless the functional dependence of the Lagrangian on $\bomega$ vanishes, i.e., 
        \begin{align} \label{eq: constraint_for_dual}
        \mathbf{b}_e^\top L(\mathbf{w})^{\dagger}-{\boldsymbol{\lambda}}_e^\top C =0,\ e\in[m].
        \end{align}
        By imposing the above condition on \eqref{eq: lagrangian function}, the dual function simplifies to $G_e(\boldsymbol{\lambda}_e)= \boldsymbol{\lambda}_e^\top \mathbf{d}$.
        Finally, by weak duality \citep{boyd2004convex}, the dual function is a tight upper bound for the primal problem, i.e., the left-hand side of \eqref{eq: max constraint}, as follows
        \begin{align}
        \underset{\bomega \in \Omega}{\mbox{max}} \ \mathbf{b}_e^\top L(\mathbf{w})^{\dagger}\bomega \leq G_e({\boldsymbol{\lambda}}_e)=\boldsymbol{\lambda}_e^\top \mathbf{d}.
        \end{align}
        The last inequality implies that the dual function $G_e({\boldsymbol{\lambda}}_e)$ can be replaced by the left-hand side of \eqref{eq: max constraint} along with the side condition \eqref{eq: constraint_for_dual}. By repeating the same procedure for the second set of constraints \eqref{eq: max constraint 2}, we will arrive at the desired equivalent problem \eqref{eq: robust weight design equivalent}. The proof is complete. $\blacksquare$
    \end{pf}
%
%\vspace{-4mm}
The robust design formulation \eqref{eq: robust weight design equivalent} is non-convex in the design variable $\mathbf{w}$, due to the pseudoinverse operation appearing in the constraints. In what follows, we propose a convex outer approximation to the feasible set of \eqref{eq: robust weight design equivalent} by virtue of the following lemma, which follows from the Schur complement condition \citep{boyd2004convex}.
\begin{lemma} \label{lem: Schur}
    Consider a weighted, connected, undirected graph $G$, with Laplacian matrix $L(\mathbf{w})$. Define $\overline{L}^{\star}$ as the minimizer of the following semidefinite program,
    \begin{align} \label{eq: Schur complement}
        \overline{L}^{\star}= & \ \underset{\overline{L}\in \mathbb{S}^{n\times n}}{\arg\min} \ \Tr(\overline{L}) \\
        \mbox{\emph{s.t. }} & \begin{bmatrix}
        L(\mathbf{w}) +\dfrac{1}{n}\mathbf{1}_n\mathbf{1}_n^\top  & I_n \\ I_n & \overline{L}\end{bmatrix} \succeq 0. \nonumber
    \end{align}
    Then, the equality $L(\mathbf{w})^{\dagger}=\overline{L}^{\star}-\frac{1}{n}\mathbf{1}_n\mathbf{1}_n^\top $ holds.
\end{lemma}
\ifx
\begin{pf}
        In Appendix \ref{appendix lem: Schur}.
\end{pf}
\fi
We now use Lemma \ref{lem: Schur} to propose a tractable convex relaxation of the optimization problem in \eqref{eq: robust weight design equivalent}. To this end, we regularize the objective function in \eqref{eq: robust weight design equivalent} with the penalty function $\Tr(\overline{L})$, and include the linear matrix inequality (LMI) in \eqref{eq: Schur complement} as an additional constraint to obtain the following outer approximation to \eqref{eq: robust weight design equivalent},
\begin{align}  \label{eq: robust weight design convex}
        \text{min}\ & f_{\mathcal{E}}(\mathbf{w})+\alpha \Tr(\overline{L}) \\ 
        \mbox{s.t. }\ &\forall e \in [m], \nonumber \\
        & \left[\begin{array}{cc}
        L(\mathbf{w})+\dfrac{1}{n}\mathbf{1}_n\mathbf{1}_n^\top & I_{n} \\
        I_{n} & \overline{L} \end{array}\right] \succeq0 , \nonumber \\
        %& f_{\mathcal{E}}\left(\mathbf{w}\right) \leq C , \nonumber \\
        %& \mbox{for} \ e\in[m] \\
            &C^\top{\boldsymbol{\lambda}}_e- \overline{L}\mathbf{b}_e=\mathbf{0}_n, \ C^\top{\boldsymbol{\gamma}}_e + \overline{L}\mathbf{b}_e =\mathbf{0}_n, \nonumber \\
            &\boldsymbol{\lambda}_e^\top \mathbf{d} \leq \sin(\gamma_d), \ \boldsymbol{\gamma}_e^\top \mathbf{d}  \leq \sin(\gamma_d), \ \boldsymbol{\lambda}_e, \  {\boldsymbol{\gamma}}_e \geq \mathbf{0}_n,\nonumber
\end{align} 
        where $\mathbf{w}\in F_{\mathbf{w}}$, $\overline{L} \in \mathbb{S}^{n\times n}$; $\boldsymbol{\lambda}_e, \boldsymbol{\gamma}_e  \in \mathbb{R}^p_{+}$ are optimization variables, and $\alpha >0$ is a regularization constant. 
        
        Notice that by Lemma \ref{lem: Schur}, for any optimal solution $(\mathbf{w}^\star,\overline{L}^\star,$ $\boldsymbol{\lambda}^\star,\boldsymbol{\gamma}^\star)$ to \eqref{eq: robust weight design convex}, the condition $\overline{L}^\star=(L(\mathbf{w}^\star)+\frac{1}{n}\mathbf{1}_n\mathbf{1}_n^\top)^{-1}$ holds if and only if $\mathbf{w}^\star$ is feasible for the original non-convex problem \eqref{eq: robust weight design equivalent}. Therefore, the regularizer coefficient $\alpha>0$ must be large enough in order to enforce the identity $\overline{L}^{\star}=(B\,\mbox{diag}(\mathbf{w}^{\star})B^{\top}+\frac{1}{n}\mathbf{1}\mathbf{1}^{\top})^{-1}$. On the other hand, too large values of $\alpha$ compromise the optimality of the objective function $f_{\mathcal{E}}(\mathbf{w})$. In practice, we use a bisection search aiming to find the smallest value of $\alpha$ for which the LMI constraint is tight.

\begin{figure*}
	\begin{center}
		\includegraphics[width=\textwidth]{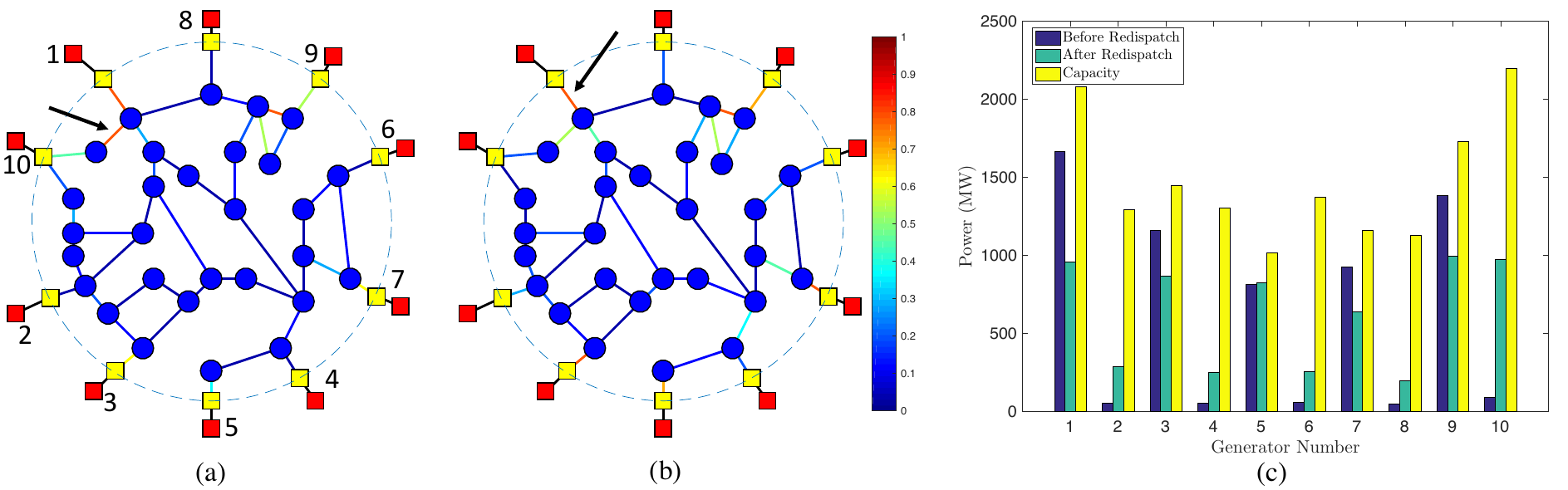}
		\caption{\small Power network diagram used in $\S$\ref{subsection:Power Redispatch}. The red squares correspond to generators; the yellow squares are terminal buses, and the blue circles are load buses. (a) Phase differences before redispatch: The most stressed link has a phase cohesiveness of $21 \deg$ (marked with a black arrow). (b) Phase differences after redispatch: The most stressed link has a phase cohesiveness of $10 \deg$ (marked with a black arrow). The colormap in each figure represents the normalized phase differences $({\theta_i^\star-\theta_j^\star})/{\varphi^\star}$ across each edge, where $\varphi^\star={\max}_{\{i,j\}\in \mathcal{E}} (\theta_i^\star-\theta_j^\star)$. (c) Distribution of power across the generators, before and  after redispatching, along with the capacity of the generators. The amount of displaced power is $53 \%$ of the total generated power.}
		\label{fig: power_redispatch_2}
	\end{center}
\end{figure*}

If there is no uncertainty in the natural frequencies $\bomega$, i.e., $\Omega=\{\bomega\}$ (which is equivalent to setting $C=[I_n, -I_n]^\top$ and $\mathbf{d}=[\bomega^\top, -\bomega^\top]^\top$ in \eqref{eq: freqeuency_uncertainty_set}), the convex relaxation of the weight design formulation in \eqref{eq: robust weight design convex} can be simplified, as stated next.
 \begin{corollary}\label{cor: robust weight design convex no uncertainty} {\normalfont \textbf{(Weight design)}} Consider the optimization problem \eqref{eq: robust weight design} with $\Omega=\{\bomega\}$ for a given $\bomega\in \mathbb{R}^n$. The following optimization problem is a convex outer approximation of \eqref{eq: robust weight design},
 	\begin{align} \label{eq: robust weight design convex no uncertainty}
 	&\underset{\mathbf{w}\in F_{\mathbf{w}},\overline{L}}{\mbox{\emph{minimize}}} \ f_{\mathcal{E}}(\mathbf{w})+\alpha \Tr(\overline{L}) \\ 
 	&\mbox{\emph{s.t.}} \left[ \begin{array}{cc}
 	L(\mathbf{w})+\dfrac{1}{n}\mathbf{1}\mathbf{1}^{\top} & I_{n} \\
 	I_{n} & \overline{L}\end{array} \right] \succeq0, \nonumber \ 
 	\|B^\top \overline{L}\bomega\|_{\infty}\leq \sin(\gamma_d). \nonumber
 	\end{align}
 \end{corollary}

\begin{remark}[Network connectivity]{
In order for Lemma \ref{lem: Schur} to be applicable in \eqref{eq: robust weight design convex}, the condition $\lambda_2(B \mbox{diag}(\mathbf{w})B^\top)>0,\ \mathbf{w}\in F_{\mathbf{w}}$ must hold. In some practical cases, we need to explicitly impose this constraint in our feasible design set (see, for example, $\S$\ref{subsection:Sparse Design}).
The following constraint can be included in the definition of $F_{\mathbf{w}}$ to guarantee a strictly positive algebraic connectivity \citep{boyd2006convex},
    \begin{align} \label{eq: connectivity constraint}
   L(\mathbf{w})+ \dfrac{\beta}{n}\mathbf{1}_n\mathbf{1}_n^\top \succeq \beta I_n, \ 0<\beta\ll1.
    \end{align}
   Notice that the eigenvalues of the left-hand side of \eqref{eq: connectivity constraint} are given by $\{\beta,\lambda_2,\ldots,\lambda_n\}$, where $\lambda_2,\ldots,\lambda_n$ are the nontrivial eigenvalues of $L(\mathbf{w})$. The above LMI enforces that $\lambda_2(L(\mathbf{w})) \geq \beta >0$; hence, the graph remains connected under \eqref{eq: connectivity constraint}. }
\end{remark}

\ifx
\subsection{Weight design with chance constraints}

In this subsection, we consider weight design problem for a network of oscillators whose natural frequencies are random with known first and second-order statistics. More precisely, consider an oscillator model with a given network structure $B$ and random natural frequencies $\bomega$ such that $\mathbb{E}[\bomega]=\bomega_0$ and $\mbox{Cov}[\bomega]=R$ are known.  Our objective is to design the link weights ($\mathbf{w}$) to achieve a desired level of phase cohesiveness $\gamma_d \in [0, \pi/2)$ with some given probability and at a minimum design cost.

\begin{problem}\textbf{\emph{(Chance-constrained weight design)}} \label{problem:}{\rm
    Assume we are given the following elements: (i) a connected undirected network with incidence matrix $B$, (ii) First and second order statistics of the natural frequencies: $\mathbb{E}[\bomega]=\bomega_0$ and $\mbox{Cov}[\bomega]=\mathbf{R}$, (iii) a convex weight-tuning cost function $f_{\mathcal{E}}\left(\mathbf{w}\right):\mathbb{R}_+^m \to \mathbb{R}$,
    (iv) a compact convex feasible design set $F_{\mathbf{w}}\subset\mathbb{R}_{+}^{m}$, (v) a desired level of phase cohesiveness $\gamma_d \in [0,\pi/2)$, and (vi) a desired maximum level of the constraint violation probability $\varepsilon \in (0,1)$. Find the optimal vector of link weights, denoted by $\mathbf{w}^{\star}$, that solves the following chance-constrained optimization problem,
    \begin{align}\label{...}
    &\mathbf{w}^{\star} =\underset{\mathbf{w}\in F_{\mathbf{w}}}{\arg\min} \ f_{\mathcal{E}}\left(\mathbf{w}\right) \\
    & \mbox{s.t. } \ \mathbb{P}\Big(\|B^\top  ( B\mbox{diag}({\mathbf{w}})B^\top )^{\dagger}\bomega \|_{\infty}\leq \sin(\gamma_d)\Big) \geq 1-\varepsilon. \nonumber
    \end{align}
    Then, by Criterion \ref{sync criterion}, the phase cohesiveness of the resulting network would satisfy $\varphi(B,\mathbf{w}^{\star},\bomega) \leq \gamma_d$ with probability of at least $1-\varepsilon$, and the minimum design cost is given by $f_{\mathcal{E}}(\mathbf{w}^{\star})$.}
\end{problem}     
In what follows, we convert the above problem into a robust weight design problem by replacing the chance constraint by a $\max$ constraint over an uncertainty set to be determined.
\begin{theorem}
Consider \eqref{eq: robust weight design} with the uncertainty set defined as 
$$\Omega=\{\bomega \in \mathbb{R}^n \colon (\bomega-\bomega_0)^\top \mathbf{R}^{-1} (\bomega-\bomega_0) \leq \varepsilon/n\}$$
Then, the phase cohesiveness of the resulting network would satisfy $\varphi(B,\mathbf{w}^{\star},\bomega) \leq \gamma_d$ with probability of at least $1-\varepsilon$, and the minimum design cost is given by $f_{\mathcal{E}}(\mathbf{w}^{\star})$.
\end{theorem}
\fi

The next section provides numerical simulations to illustrate the effectiveness of the proposed framework in designing optimal networks of oscillators from the point of view of phase cohesiveness.

%!TEX root = main_Automatica.tex
\section{Applications}
\label{section:NUMERICAL SIMULATIONS}

This section illustrates the use of our optimization framework in several problems of practical interest, namely, power re-dispatch in electric grids ($\S$\ref{subsection:Power Redispatch}), sparsity-promoting network design ($\S$\ref{subsection:Sparse Design}), robust network design for distributed analog clocks ($\S$\ref{subsection: robust design}), and the Braess' paradox ($\S$\ref{example: Braess's Paradox}).

%{\color{blue} We start by considering a power redispatch problem in elecric grids as an application of our solution to the frequency design problem (Problem \ref{prob:Optimal Budget Problem-Frequency design}). Next, we study the problem of finding the sparsest network able to achieve a desired level of phase cohesiveness. This problem is an application of weight design problem (Problem \ref{problem: Budget-Constrained Weight Design}). We then study the design of networks to synchronize distributed analog clocks with uncertain frequencies. Lastly, we use our framework to study the Braess' paradox in the context of power systems using our solution to the weight design problem.}

\begin{figure*}[t]
	\begin{center}
		\includegraphics[width=\textwidth]{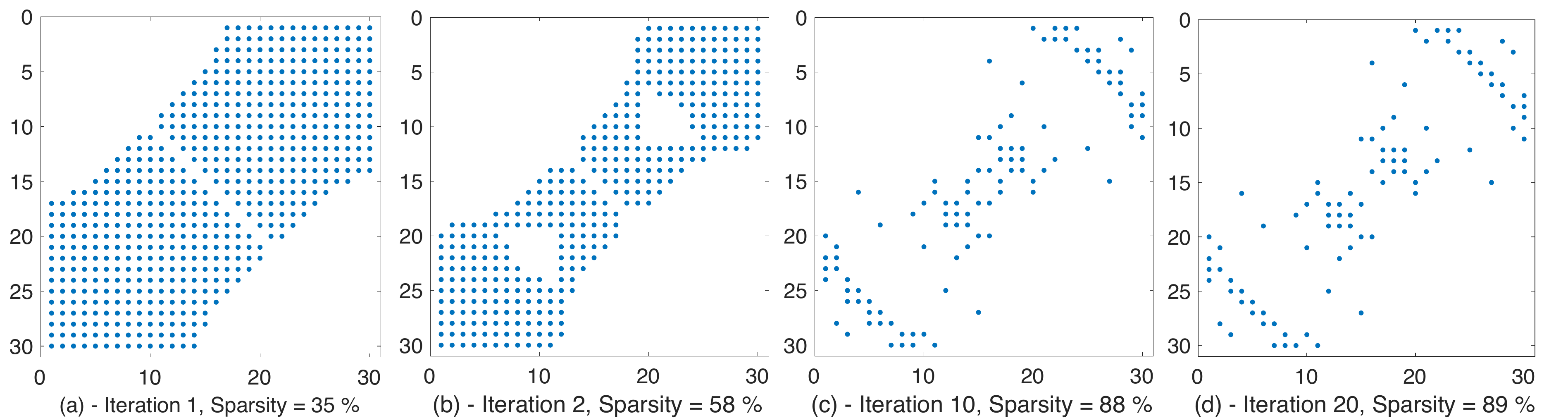}
		\caption{Evolution of the sparsity pattern for the numerical example in $\S$\ref{subsection:Sparse Design}. In (a), we plot the sparsity pattern of the adjacency matrix after the first iteration, which contains $65 \%$ of the candidate edges (or, equivalently, presents a 35\% of sparsity). Figures (b), (c) and (d) represent the networks obtained at iterations $2$, $10$, and $20$, respectively.}
		\label{fig: Sparse_Network_Design_1}
	\end{center}     
\end{figure*}

\subsection{Power redispatch/load shedding}\label{subsection:Power Redispatch}

Power redispatch (respectively, load shedding) refers to the process of adjusting the generators power injection (respectively, the loads power consumption) to relieve overloaded transmission lines, or to re-balance the grid after a fault or unexpected event. In this situation, the redistribution of power injections (or consumptions) is commonly used as a short-term remedial action to resolve congestions and balancing issues.

A lossless power network is typically modeled via \eqref{eq:Kuramoto model}, where the set of generator nodes are denoted by $\mathcal{V}_1$ and the set of load buses are denoted by $\mathcal{V}_2$. The steady-state operating point satisfies \eqref{eq: fixed point}, where $\omega_i \in \mathbb{R}$ is the net power injected into node $i \in [n]$. The incidence matrix $B$ represents the connectivity of the network, and the $e$-th edge weight is $w_e=|V_i||V_k|Y_{ik}$, where $V_i=|V_i| \exp(\mathbf{j} \theta_i)$ and $V_k=|V_k| \exp(\mathbf{j} \theta_k)$  are complex voltages at nodes $i$ and $k$, $Y_{ik}>0$ is the susceptance of the transmission line $\{i,k\} \in \mathcal{E}$, and $\mathbf{j}$ denotes the unit imaginary number.

Assume that the vector of net power allocations $\bomega_0 \in \mathbb{R}^n$ (with $\mathbf{1}_n^\top \bomega_0=0$) of an electric grid is such that the phase cohesiveness does not satisfy a desired threshold $\gamma_d \in [0,\pi/2)$. In this situation, we are interested in modifying the power allocation by a vector of increments $\Delta \bomega \in \mathbb{R}^n$ in order to satisfy the phase cohesiveness level $\gamma_d$ while minimizing the total redispatch/load shedding cost. In mathematical terms, we aim at solving the following optimization problem (see \eqref{eq: Cohesiveness-Constrained Frequency Design}):
\begin{align} \label{eq: power_redispatch} 
%\begin{aligned} 
\underset{\Delta \bomega}{\min} &\sum_{k \in \mathcal{V}} g_k(\Delta \omega_k) \nonumber \\
\mbox{s.t. }
&\|B^\top L(\mathbf{w})^\dagger (\bomega_0+\Delta \bomega)\|_{\infty}\leq \sin(\gamma_d), \nonumber \\ 
&\mathbf{1}_n^\top  \Delta \bomega=0, \ \underline{\bomega} \leq \bomega_0 + \Delta \bomega \leq \overline{\bomega},
%\end{aligned}
\end{align}
where $\underline{\bomega} \leq \overline{\bomega} \in \mathbb{R}^n$ are the vectors of admissible lower and upper bounds on the nodal power injections. The convex objective function $g_k(\Delta \omega_k)$ is the redispatch/load shedding cost at node $k \in \mathcal{V}$. %The first constraint imposes the desired level of phase cohesiveness, the second constraint guarantees the power balance across the grid, and the last constraint prescribes limits on the capacity of generators and loads.  
In the optimal redispatch problem, we are allowed to adjust the power injected in the generators buses only, such that $\Delta \omega_k=0$ for all $k \in \mathcal{V}_2$. In contrast, in the load shedding problem, we can adjust the load buses only, such that $\Delta \omega_k=0$ for all $k \in \mathcal{V}_1$.

In our numerical evaluation, we consider a power redispatch problem for the New England power grid depicted in Fig. \ref{fig: power_redispatch_2}-(a) \citep{dorfler2010spectral}. The network data $(B,\mathbf{w},\{\omega_{0,i}\}_{i \in \mathcal{V}_2},\underline{\bomega},\overline{\bomega})$ are obtained from \cite{zimmerman2011matpower}, assuming that the transmission lines are lossless and there are no transformers or phase shifters. There are 10 generators (red squares) connected by terminal buses (yellow squares) to 29 load buses (blue circles). We assume that the generators $\{1,3,5,7,9\}$ are generating power at $95 \%$ of their capacity, while the generators $\{2,4,6,8,10\}$ are generating at $5 \%$ of their capacity. For these particular values, we numerically solve the fixed point equation \eqref{eq: fixed point}, from where we obtain a value of phase cohesiveness of $21 \deg$. We then solve the problem of redispatching the minimum amount of power to guarantee a phase cohesiveness of $\gamma_d=10 \deg$ by solving \eqref{eq: power_redispatch} with $g_k(\Delta \omega_k)=|\Delta \omega_k|$ and $\Delta \omega_k=0$ for $k \in \mathcal{V}_2$. The resulting power flow distribution is depicted in Fig. \ref{fig: power_redispatch_2}-(b). The total redispatched power $\|\Delta \bomega \|_1$ is $53 \%$ of the total power generation in the network.

%\begin{figure}
%	\begin{center}
%		\includegraphics[width=0.47\textwidth]{power_dispatch_2.pdf}
%		\caption{\small Power network diagram used in $\S$ \ref{subsection:Power Redispatch}. The red squares correspond to generators; the yellow squares are terminal buses, and the blue circles are load buses. (a) Before redispatch: The most stressed link has a phase cohesiveness of $21 \deg$. (b) After redispatch: The most stressed link has a phase cohesiveness of $10 \deg$. The colormap in each figure represents the normalized phase differences $({\theta_i^\star-\theta_j^\star})/{\varphi^\star}$ across the edges, where $\varphi^\star=\underset{\{i,j\}\in \mathcal{E}}{\max} (\theta_i^\star-\theta_j^\star)$. }
%		\label{fig: power_redispatch}
%	\end{center}
%\end{figure}

%\subsection{Sparse network design}
%Consider the model \eqref{eq: Kuramoto model-matrix} for a connected undirected graph $G=(\mathcal{V},\mathcal{E}_0)$ with given incidence matrix $B_0$ and vector of natural frequencies $\bomega_0$. We consider the following design problem: Given the boolean structure of the network, tune the link weights, within a closed convex set $F_{\mathbf{w}} \subseteq \mathbb{R}_{+}^m$, such that the sum of the weights is minimized and the resulting network is sychronizable with a specified phase cohesiveness level $\gamma_d \in [0,\pi/2)$. In other words, the goal is to solve the following problem,
%
%
\subsection{Sparse network design}\label{subsection:Sparse Design}

Consider the network dynamics in \eqref{eq: Kuramoto model-matrix} for a connected undirected graph $G=(\mathcal{V},\mathcal{E})$ with given $B$ and $\bomega$. We consider the weight design problem \eqref{eq: robust weight design} (and its convex approximation \eqref{eq: robust weight design convex no uncertainty}) where the cost function is given by the sum of the edge weights, i.e., $f_{\mathcal{E}}(\mathbf{w})=\|\mathbf{w}\|_1$. Because of its sparsity-promoting nature, the $\ell_1$ norm is typically used in sparse design problems, where a solution with many zero entries is desired. In practice, however, the designed network might have a relatively large number of links with optimal weights close to zero, but not exactly zero. To promote sparsity, we propose to use the \emph{re-weighted} $\ell_1$ minimization algorithm, described in \cite{candes2008enhancing}. In this algorithm, a sequence of weighted $\ell_1$-norm problems are solved such that, in each round, the weights of the $\ell_1$ norm are updated to promote sparsity in the next round. The re-weighted $\ell_1$ minimization algorithm is summarized in Algorithm \ref{tab:weight design for optimal cohesiveness}, and described below. 

\emph{Description of the algorithm}: In Step 1, the coefficients of the $\ell_1$ norm are initialized at one (i.e., $\mathbf{p}^{(1)}=\mathbf{1}_{m}$), and the incidence matrix of the selected (nonzero) edges is set to ${B}^{(1)}_s=B$. In Step 3 of iteration $k$, the convex semidefinite relaxation \eqref{eq: weighted L1 weight design} is solved in order to obtain the optimal edge weights $\mathbf{w}^{(k)}$. In Step 4, the components of $\mathbf{p}^{(k)}$ are updated inversely proportional to the corresponding components of $\mathbf{w}^{(k)}$. The constant $0<\varepsilon \ll 1$ is used to avoid singularities. In Step 5, ${B}_s$ is updated to include only the selected (nonzero) edges obtained at Step 3 for the next iteration. Steps 2 to 6 are repeated for a specified number of iterations (denoted by $k_{\max}$) or until a desired sparsity is achieved. The incidence matrix $B_s^{k_{\max}}$ will then include the final selected edges.

\begin{figure*}[t]
	\centering
	\includegraphics[width=0.98\textwidth]{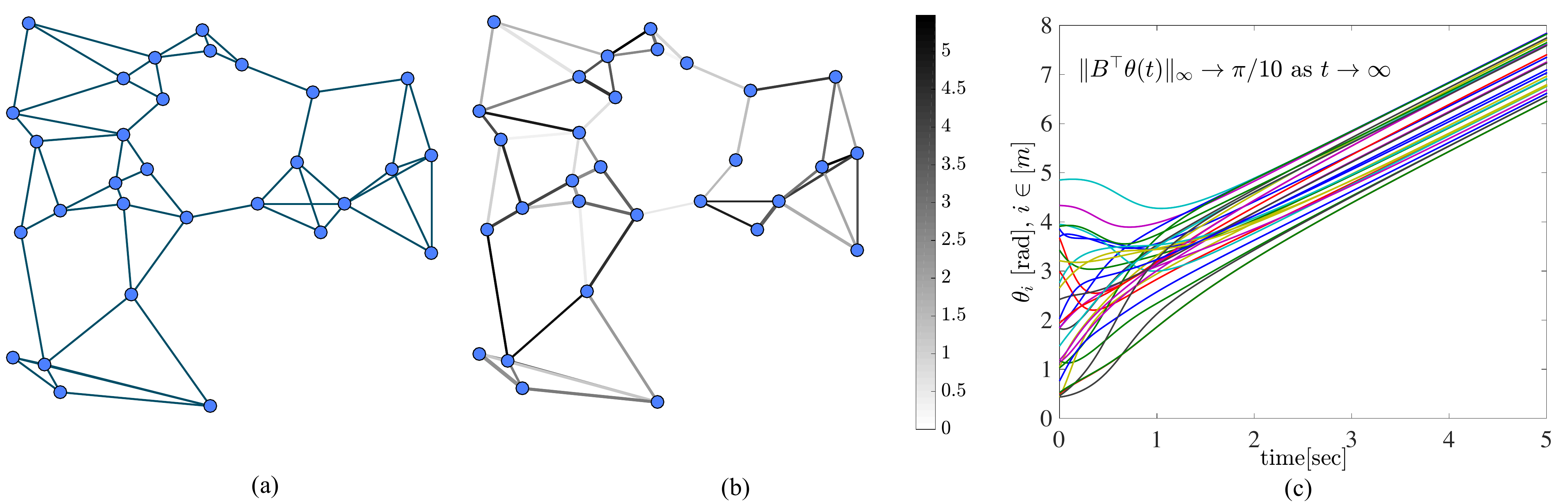}
	\caption{Sensor network considered in $\S$\ref{subsection: robust design}. The initial graph is plotted in (a), the designed graph in (b), and the time evolution of phases for the worst-case realization of the natural frequencies in (c).}
	\label{fig: sensor network}
\end{figure*}

\begin{algorithm}
	{\footnotesize \vspace*{1ex}
		\textbf{Given}: $B=[\mathbf{b}_{e}]_{\in \mathcal{E}}$, $\bomega$, $F_{\mathbf{w}}$, $ \gamma_d$, $\alpha$, $k_{\max}$, and $0<\varepsilon \ll 1$.\\
		\begin{algorithmic}[1]
			\State set $\mathbf{p}^{(1)}=\mathbf{1}_{m}$ and ${B}^{(1)}_s=B$; \label{step 1}
			\For {$k=1,\ldots,k_{\max}$} 
			\State solve \eqref{eq: weighted L1 weight design} to obtain $\mathbf{w}^{(k)}$:
			\begin{align} \label{eq: weighted L1 weight design}
			\mathbf{w}^{(k)}=& \arg \underset{\mathbf{w}\in F_{\mathbf{w}},\overline{L}}{\mbox{min}} \ \sum_{e=1}^{m} p_e^{(k)}|w_{e}|+\alpha \Tr(\overline{L}) \\ 
			\mbox{{s.t.}} & \left[ \begin{array}{cc} 
			B\mbox{{diag}}(\mathbf{w})B^{\top}+\dfrac{1}{n}\mathbf{1}\mathbf{1}^{\top} & I_{n} \\
			I_{n} & \overline{L}\end{array} \right] \succeq0, \nonumber \\
			& \|{B_s^{(k)}}^\top \overline{L} \bomega\|_{\infty} \leq \sin(\gamma_d), \nonumber
			%& \|{B}_s^{(k)\top} \overline{L}\bomega_0\|_{\infty}\leq \sin(\gamma_d). \nonumber
			\end{align}
			\State update $p_{e}^{(k+1)}=(\varepsilon+w_{e}^{(k)})^{-1},\ e \in [m]$;
			\State update ${B}^{(k+1)}_s=\left[\mathbf{b}_{e}\right]_{\{e \in \mathcal{E} \colon w_{e}^{(k)}>0 \}}$;
			%\STATE update $\mathbf{w}=\mathbf{w}^{+}$.
			%\STATE update $c_i=(\varepsilon+w_i)^{-1},\ i \in [m]$.
			%\STATE{update $\mathbf{w}^{+}$ as the optimal solution of \eqref{eq: weighted L1 weight design}}.
			\EndFor
		\end{algorithmic}}
		\caption{\hspace*{-.5ex}: sparse weight design} \label{tab:weight design for optimal cohesiveness}
	\end{algorithm}

\ifx	
	\begin{figure}
		\begin{center}
			\includegraphics[width=0.35\textwidth]{Sparse_Network_Design_2.pdf}
			\caption{Time evolution of phases  for the numerical example in $\S$ \ref{subsection:Sparse Design}. The phase cohesiveness satisfies the design requirement, i.e., $\varphi(B,\mathbf{w}^\star,\bomega_0)=30 \deg$.}
			\label{fig:Sparse_Network_Design_2}
		\end{center}
	\end{figure}
\fi
	
	In our numerical experiments, we assume that $n=30$, $B=B_{K_n}$ where $K_n$ denotes the all-to-all graph, $\omega_{i}=-1+2\frac{i-1}{n-1}$ for $\ i\in [n]$, $m=\binom{n}{2}$, $F_{\mathbf{w}}=\mathbb{R}_{+}^m$, $\gamma_d=30\deg$, and $\alpha=0.5$. In other words, the network designer is allowed to connect any pair of nodes. To maintain the connectivity of the network, we include the LMI in \eqref{eq: connectivity constraint} with $\beta=10^{-4}$ in the definition of $F_{\mathbf{w}}$. Fig. \ref{fig: Sparse_Network_Design_1} illustrates the evolution of the sparsity pattern of the adjacency matrix as Algorithm \ref{tab:weight design for optimal cohesiveness} progresses.
	
%\ifx
	\subsection{Robust synchronization of distributed analog clocks} \label{subsection: robust design}
	Consider a wireless sensor network consisting of $n$ processors $\mathcal{V}=[n]$ equipped with analog clocks. In order to efficiently perform distributed computations across the network, the clocks are required to synchronize their phases. The oscillator model \eqref{eq:Kuramoto model} without inertia (i.e., $\mathcal{V}_1=\emptyset,\ d_i=0,\ i\in \mathcal{V}_2$) can be used as a distributed synchronization scheme for synchronizing the phases \citep{simeone2008distributed}. In this context, the matrix $B \in \mathbb{R}^{n\times m}$ is the incidence matrix of the communication graph, $\mathbf{w} \in \mathbb{R}^{m}_{+}$ is the vector of connection strengths, and $\bomega \in \mathbb{R}^n$ is the vector of natural frequencies of the clocks. In practice, the natural frequencies $\omega_i, \ i\in[n]$ are uncertain due to hardware imperfections and aging. Therefore, the communication graph must be designed in order to synchronize the clocks in the presence of uncertainties in the natural frequencies. More specifically, we aim to allocate the minimum amount of edge weights while guaranteeing a desired level of phase cohesiveness. We pose this allocation problem as Problem \ref{problem: Budget-Constrained Weight Design} with
	\ifx
	\begin{align}
	\mathbf{w}^\star(\gamma_d) =\arg &\underset{\mathbf{w}\in F_{\mathbf{w}}}{\min} \|\mathbf{w}\|_{1} \\
	& \mbox{\:  s.t.  \:} \max_{\bomega \in \Omega}\varphi(B,\mathbf{w},\bomega) \leq \gamma_d, \nonumber
	\end{align}
	\fi 
	the sum of the weights $\sum_{e=1}^{m} w_e=\|\mathbf{w}\|_1$ as the cost function. In our numerical simulations, we consider the sensor network depicted in Fig. \ref{fig: sensor network}-(a)  with $n=30$ processors and $m=56$ links. We assume that the natural frequencies are nominally equal to $1$ with $20\%$ uncertainty, i.e., 
	\begin{align} \label{eq: freqeuency_uncertainty_set_1}
	\Omega=\{\bomega \in \mathbb{R}^n \colon \  0.8 \, \mathbf{1}_n \leq \bomega \leq 1.2 \, \mathbf{1}_n\}.
	\end{align}
	This box constraint set can be written in the polyhedral form \eqref{eq: freqeuency_uncertainty_set} with $C=[I_n, -I_n]^\top$ and $\mathbf{d}=[1.2\mathbf{1}_n^\top, -0.8 \mathbf{1}_n^\top]^\top$. We then solve \eqref{eq: robust weight design convex} with $\gamma_d=\pi/10 \ \mbox{rad}$, and the feasible design set being the positive orthant, $F_{\mathbf{w}}=\mathbb{R}_{+}^{m}$.
	The resulting network is illustrated in Fig. \ref{fig: sensor network}-(b), where the optimal cost is $\|\mathbf{w}^\star\|_1 \approx 70$. To verify the robustness of the designed network, we numerically integrate \eqref{eq:Kuramoto model}, using a realization from the worst-case set of natural frequencies $\Omega^\star \subset \Omega$, defined as
	$\Omega^\star = \argmax_{\bomega \in \Omega} \|B^\top (B\mbox{diag}(\mathbf{w^{\star}})B^\top)^{\dagger} \bomega\|_{\infty}$. 
	The resulting evolution is plotted in Fig. \ref{fig: sensor network}-(c). We observe that $\lim_{t\to\infty}\|B^\top \theta(t)\|_{\infty}=\pi/10$, as expected; hence, the phase cohesiveness of the optimal network is guaranteed to be less that $\pi/10$ for all $\bomega \in \Omega$.
	%
	%\begin{figure} 
	%    \centering
	%    \includegraphics[width=0.4\textwidth]{Robust_Weight_Design_2.pdf}
	%    \caption{Sensor network}
	%    \label{fig: sensor network phases}
	%\end{figure}
%\fi	

	\subsection{Braess' paradox in power systems} \label{example: Braess's Paradox}
\begin{figure}
		\begin{center}
			\includegraphics[width=0.47\textwidth]{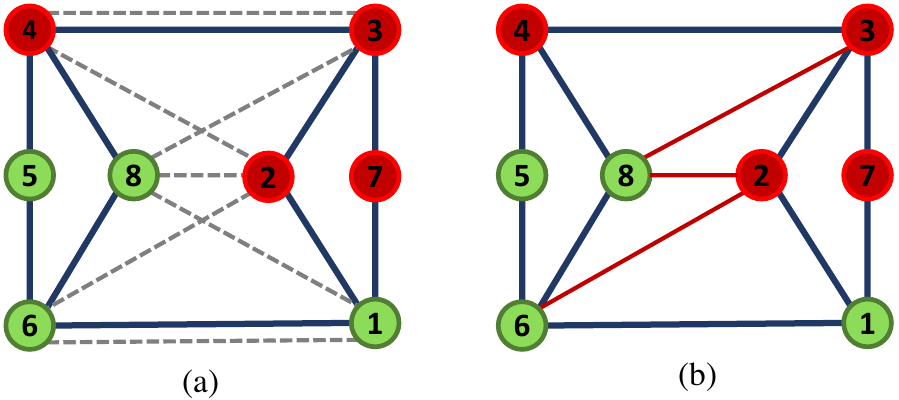}
			\caption{\small Power network considered in $\S$\ref{example: Braess's Paradox} (from \cite{witthaut2012braess}). Load buses and generators are depicted as green and red nodes, respectively. Dashed lines in (a) represent candidate edges that can be added to an existing network, while solid lines are the already existing edges. In (b), the three red lines denote the lines added after the optimization.}
			\label{fig:Braess_Paradox_1}
		\end{center}
	\end{figure}

	The Braess' paradox refers to the counter-intuitive phenomenon of losing synchrony as a result of adding new links to a network, or strengthening the existing ones \citep{witthaut2012braess}. To illustrate this paradox, we consider the lossless power network represented in Fig. \ref{fig:Braess_Paradox_1}-(a), which we will refer to as $G_0$. This network has $4$ generators (orange nodes), $4$ load buses (green nodes), and $m_0=10$ transmission lines (solid lines). All nodes are assumed to have the same value of power demand/generation, in particular,  $\omega_{i}=0.95$ for generators and $\omega_{i}=-0.95$ for load buses. Furthermore, all the edges in $G_0$ are assumed to have identical capacity equal to $1$. For these numerical values, the phase cohesiveness satisfies $\sin(\varphi(B_0,\mathbf{w}_0,\bomega))=0.95$, where $B_0$ is the incidence matrix of $G_0$ and $\mathbf{w}_0=\mathbf{1}_{10}$.  
	For this network, let us consider the problem of adding new lines to the network (chosen from a set of candidate edges) in order to decrease the phase cohesiveness below the value $\gamma_d=\pi/3$. The candidate lines are indicated by dashed lines in Fig. \ref{fig:Braess_Paradox_1}-(a). We denote the subgraph induced by the candidate lines as $G_c$, its incidence matrix as $B_c \in \mathbb{R}^{8 \times 7}$, and its weights as $\mathbf{w}_c\in\mathbb{R}_{+}^7$. In what follows, we minimize the total capacity (measured as the $\ell_1$ norm of $\mathbf{w}_c$) added to the network, which can be posed as the following optimization problem:
	\begin{align}\label{eq: l1_norm_Braess}
	\mathbf{w}_c^{\star} =&\arg \underset{\mathbf{w}_c\in F_{\mathbf{w}_c}}{\min} \|\mathbf{w}_c\|_{1} \\
	%& \mbox{ \; s.t. } \; \; \varphi(B_s,\mathbf{w}_s,\bomega_0) \leq \gamma_d, \nonumber
	\mbox{ \; s.t. } \; \; &\|B_0^\top L^\dagger \bomega\|_{\infty} \leq \sin(\gamma_d),
	\|B_c^\top L^\dagger \bomega\|_{\infty} \leq \sin(\gamma_d), \nonumber \\
	& L=B_0\mbox{diag}(\mathbf{w}_0)B_0^\top+B_c\mbox{diag}(\mathbf{w}_c)B_c^\top. \nonumber
	\end{align}
	By \eqref{eq: robust weight design convex no uncertainty}, the corresponding relaxation is
	\begin{align} \label{eq: l1_norm_Braess_SDP}
	\mathbf{w}_c=& \arg \underset{\mathbf{w}_c\in F_{\mathbf{w}},\overline{L}}{\mbox{min}} \ \|\mathbf{w}_c\|_1+\alpha \Tr(\overline{L}) \\ 
	\mbox{{s.t.}} &\left[ \begin{array}{cc} 
	L_0+B_c\mbox{diag}(\mathbf{w}_c)B_c^\top+\dfrac{1}{n}\mathbf{1}\mathbf{1}^{\top} & I_{n} \\
	I_{n} & \overline{L}\end{array} \right] \succeq0, \nonumber \\
	& \|B_0^\top \overline{L}\bomega\|_{\infty}\leq \sin(\gamma_d), \|{B}_c^\top \overline{L}\bomega\|_{\infty}\leq \sin(\gamma_d), \nonumber
	\end{align}
	where $L_0=B_0\mbox{diag}(\mathbf{w}_0)B_0^\top$. The resulting network is depicted in Fig \ref{fig:Braess_Paradox_1}-(b). The optimal nonzero edges are $w_{28}$, $w_{38}$, and $w_{26}$, and the remaining candidate links ($w_{34}$, $w_{24}$, $w_{16}$, and $w_{18}$) have zero optimal value. 
	
\begin{figure}
		\begin{center}
			\includegraphics[width=0.45\textwidth]{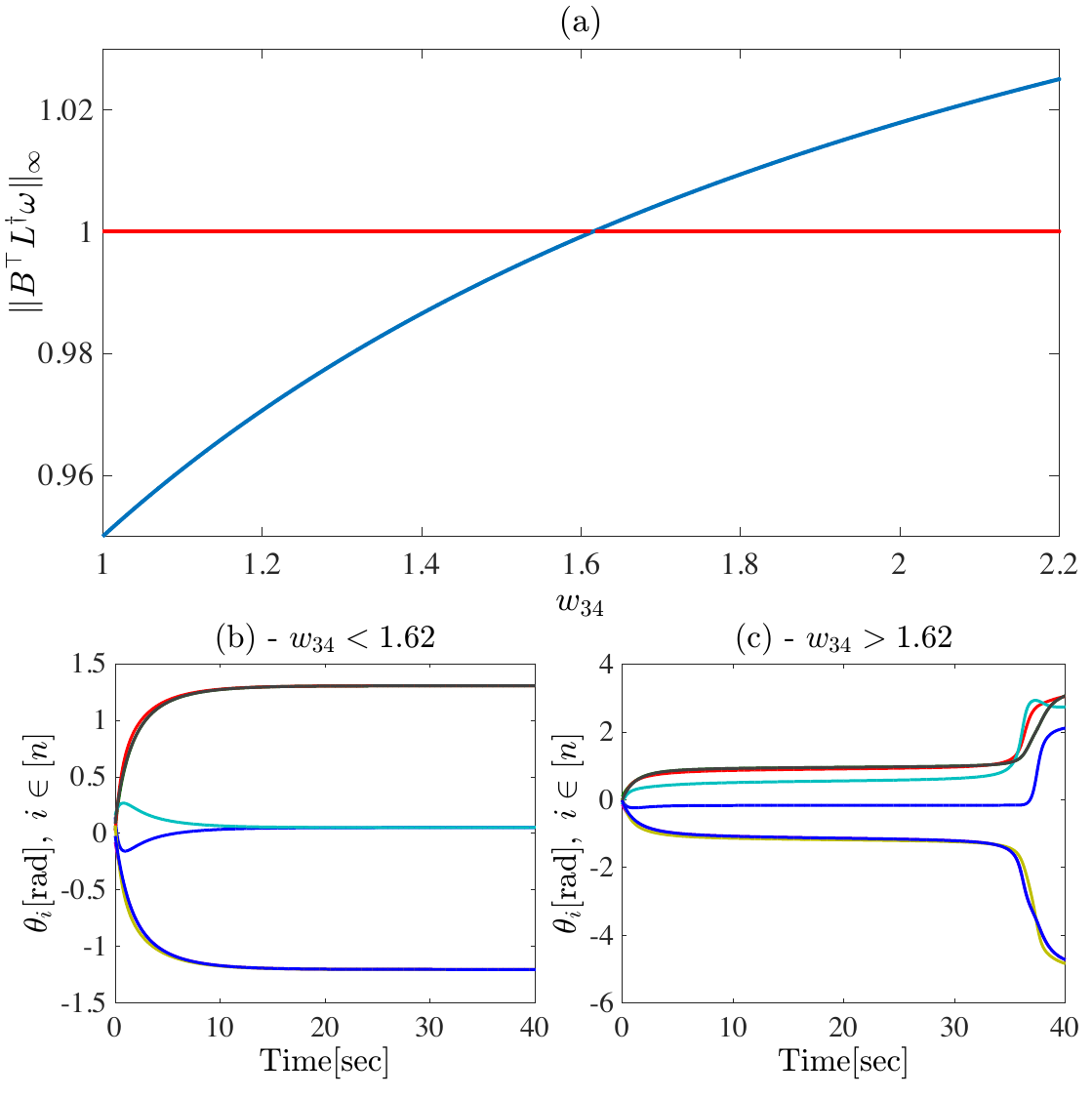}
			\caption{(a) Variation of $\|B^\top L^{\dagger} \bomega \|_{\infty}$ as a function of $w_{34}$. (b) Time evolution of phases when $w_{34}<1.62$. (c) Time evolution of phases when $w_{34}>1.62$.}
			\label{fig:Braess_Paradox_2}
		\end{center}
	\end{figure}

	To relate the result of our optimization to the Braess' paradox, we run the following experiment: increase the capacity $w_{34}$ and plot the variation of $\|B_0^\top L^\dagger \bomega \|_{\infty}$ as a function of $w_{34}$. This variation is plotted in Fig \ref{fig:Braess_Paradox_2}-(a), where we observe how, as we increase the link strength $w_{34}$, the value of $\|B_0^\top L^\dagger \bomega \|_{\infty}$ increases monotonically and crosses the stability threshold at $w_{34}\approx 1.62$. In other words, increasing the value of $w_{34}$ has a detrimental effect on the network stability. To validate our claims, we plot the time evolution of the phase dynamics for $w_{34}=1<1.62$ (Fig. \ref{fig:Braess_Paradox_2}-(b)) and $w_{34}=2>1.62$ (Fig. \ref{fig:Braess_Paradox_2}-(c)), in which we observe how the network dynamics transition from a stable to an unstable regime. Similar results can be observed when we increase $w_{16}$ or add the new lines $w_{24}$ and $w_{28}$.
	\ifx
	Similarly, we consider the effect of adding a new line connecting nodes 2 and 4 by increasing its capacity from $0$ to $1$. In Fig. \ref{fig:Braess_Paradox_3}, we plot the variation of $\|B_0^\top L^\dagger \bomega \|_{\infty}$ as we increase  $w_{24}$, and observe that when $w_{24}\approx 0.41$ the network transition from a stable to an unstable state. We validate this claim by plotting the time evolution of the phase dynamics for $w_{24}=0<0.41$ (Fig. \ref{fig:Braess_Paradox_3}-(b)) and $w_{34}=0.8>0.41$ (Fig. \ref{fig:Braess_Paradox_3}-(c)).  Similar results can be obtained when we increase $w_{28}$.
	\fi
	These observations confirm that the proposed optimization problem \eqref{eq: l1_norm_Braess_SDP} has assigned zero weight to those links that are detrimental to the phase cohesiveness. More generally, our optimization framework, which is based on Criterion \ref{sync criterion}, is capable of identifying those lines inducing the Braess' paradox.
	\ifx
	\begin{figure}
		\begin{center}
			\includegraphics[width=0.47\textwidth]{Braess_Paradox_3.pdf}
			\caption{(a) Variation of $\|B^\top L^{\dagger} \bomega \|_{\infty}$ as a function of $w_{24}$. (b) Time evolution of phases when $w_{24}<0.41$. (c) Time evolution of phases when $w_{24}>0.41$.}
			\label{fig:Braess_Paradox_3}
		\end{center}
	\end{figure}
	\fi
	%
	%
	%\begin{figure}
	%    \centering
	%    \begin{subfigure}[b]{\linewidth}
	%        \includegraphics[width=0.5\textwidth]{Braess_Paradox_1_1.pdf}
	%        \caption{}
	%    \end{subfigure}\par\vfill \bigskip
	%    \begin{subfigure}[b]{\linewidth}
	%        \includegraphics[width=0.5\linewidth]{Braess_Paradox_1_1.pdf}
	%        \caption{} \label{fig_two_agent_error}
	%    \end{subfigure}\bigskip
	%    \caption{}
	%\end{figure}
	%
	
	%\begin{figure}[H]
	%        \centering
	%        \includegraphics[width=0.5\textwidth]{Braess_Paradox_3.pdf}
	%        \caption{Variation of phase cohesiveness for the network of Example \ref{example: Braess's Paradox} as a new link between node 2 and 4 is added.}
	%        \label{fig:Braess_Paradox_3}
	%\end{figure}

	%\begin{figure}
	%        \centering
	%        \includegraphics[width=0.5\textwidth]{Braess_Paradox_4.pdf}
	%        \caption{Temporal evolution of phase angles for the network of Example \ref{example: Braess's Paradox} when $w_{34}<1$ (top left), $w_{34}=1$ (top right), $w_{34}>1$ (bottom left), and  $w=w^{\star}$ (bottom right).}
	%        \label{fig:Braess_Paradox_4}
	%\end{figure}

%!TEX root = main_Automatica.tex
\section{Conclusions}
\label{section:CONCLUSIONS}
This paper proposes a convex optimization framework for designing the natural frequencies and the coupling weights in a network of nonidentical coupled oscillators. We have used phase cohesiveness as our design constraint, capturing both the steady-state performance and the stability of the network. 
In this context, we have addressed the following network design problems:
(\emph{i}) the nodal-frequency design problem, in which we design the natural frequencies of the oscillators for a given network, and 
(\emph{ii}) the edge-weight design problem, in which we design the edge weights.
For the latter case, we have also developed a robust framework to design networks under frequency uncertainty, in which the uncertainty model is deterministic and set-based. We have illustrated the applicability of our results using several network design problems of practical interest, namely, a power redispatch case study in power grids ($\S$\ref{subsection:Power Redispatch}), a sparsity-promoting design problem ($\S$\ref{subsection:Sparse Design}), a robust network design problem in the context of distributed analog clocks ($\S$\ref{subsection: robust design}), and a network design problem in which we illustrate the Braess' paradox ($\S$\ref{example: Braess's Paradox}).
%\appendix
%\input{Appendix.tex}


%%%%%%%%%%%%%%%%%%%%%%%%%%%%%%%%%%%%%%%%%%%%%%%%%%%%%%%%%%%%%%%%%%%%%%%%%%%%%%%%

%%%%%%%%%%%%%%%%%%%%%%%%%%%%%%%%%%%%%%%%%%%%%%%%%%%%%%%%%%%%%%%%%%%%%%%%%%%%%%%%

%%%%%%%%%%%%%%%%%%%%%%%%%%%%%%%%%%%%%%%%%%%%%%%%%%%%%%%%%%%%%%%%%%%%%%%%%%%%%%%%

%%%%%%%%%%%%%%%%%%%%%%%%%%%%%%%%%%%%%%%%%%%%%%%%%%%%%%%%%%%%%%%%%%%%%%%%%%%%%%%%

%\bibliographystyle{harvard}
\bibliography{references}

\end{document}